\documentclass[10pt]{amsart}

\usepackage{amssymb,epsfig}
\usepackage{amsmath,verbatim}
\newcommand{\nc}{\newcommand}
\nc{\nt}{\newtheorem}
\nc{\bs}{\bigskip}
\nc{\dmo}{\DeclareMathOperator}

\nt{main}{Theorem}
\nt{thm}{Theorem}[section]
\nt{prop}[thm]{Proposition}
\nt{lem}[thm]{Lemma}
\nt{fact}[thm]{Fact}
\nt{cor}[thm]{Corollary}
\nt{conj}[thm]{Conjecture}
\nt{question}[thm]{Question}
\nt{problem}[thm]{Problem}
\nt{remark}[thm]{Remark}
\nt{result}[thm]{}

\nc \A{\mathcal A}
\nc{\C}{\mathcal{C}}
\nc {\RR}{\mathcal R}
\nc{\LL}{\mathcal{L}}
\nc{\SL}{\mathcal{SL}}
\nc {\W}{{\mathcal W}}
\nc {\UR}{\mathcal {UR}}
\nc {\VR}{\mathcal {VR}}
\nc {\UW}{\mathcal {UW}}
\nc {\WR}{\mathcal {WR}}
\nc {\AL}{\mathcal {AL}}
\nc {\SW}{\mathcal {SW}}
\nc {\SA}{\mathcal {SA}}
\nc {\SR}{\mathcal {SR}}
\nc{\PR}{\mathcal{PR}}
\nc{\PUR}{\mathcal{PUR}}
\nc{\PW}{\mathcal{PW}}
\nc{\PUW}{\mathcal{PUW}}

\nc{\ur}{U\!R}
\nc{\vr}{V\!R}
\nc{\pur}{PU\!R}

\nc{\Z}{\mathbb{Z}}
\nc{\R}{\mathbb{R}}
\nc{ \smZ}{\mathbb Z}
\nc {\Diff}{\hbox{\rm Diff}}
\nc {\rel}{\hbox{\ \rm rel}\,} 
\nc {\cln}{\,\colon\!} % Used for function notation.

\nc \RP{\R\hbox{\rm P}}

\nc \e{\varepsilon}

\nc \bdy{\partial}

\nc \incl{\hookrightarrow}
\nc \To{\longrightarrow}

\nc{\p}[1]{\medskip\paragraph{{\bf #1}}}
\nc{\margin}[1]{\marginpar{\scriptsize #1}}

\title{Configuration Spaces of Rings and Wickets}

\include{diagram}

\begin{document}
	
\input{epsf.sty}

\author{Tara E. Brendle}

\author{Allen Hatcher}

\address{Tara E. Brendle \\ Department of Mathematics\\University of Glasgow\\
University Gardens\\Glasgow G12 8QW  UK\\ t.brendle@maths.gla.ac.uk}

\address{Allen Hatcher \\ Department of Mathematics \\ Malott Hall \\ Cornell University \\ Ithaca, NY  14853-4201  USA \\ hatcher@math.cornell.edu}

\thanks{The first author gratefully acknowledges support from the National Science Foundation.}

%\email{t.brendle@maths.gla.ac.uk}

%\email{hatcher@math.cornell.edu}



\begin{abstract}
The main result in this paper is that the space of all smooth links in $\R^3$ isotopic to the trivial link of $n$ components has the same homotopy type as its finite-dimensional subspace consisting of configurations of $n$ unlinked Euclidean circles (the `rings' in the title). There is also an analogous result for spaces of arcs in upper half-space, with circles replaced by semicircles (the `wickets' in the title). A key part of the proofs is a procedure for greatly reducing the complexity of tangled configurations of rings and wickets. This leads to simple methods for computing presentations for the fundamental groups of these spaces of rings and wickets as well as various interesting subspaces. The wicket spaces are also shown to be aspherical.
\end{abstract}

\maketitle

\section{Introduction}

The classical braid group $B_n$ can be defined as the fundamental group of the space of all configurations of $n$ distinct points in $\R^2$. In this paper we consider a $3$-dimensional analog which we call the {\it ring group\/} $R_n$.  This is the fundamental group of the space $\RR_n$ of all configurations of $n$ disjoint pairwise unlinked circles, or rings, in $\R^3$, where we mean the word `circle' in the strict Euclidean sense.  It is not immediately apparent that $\RR_n$ is path-connected, but in Section~2 we recall a simple geometric argument from [FS] that proves this. In particular, this says that configurations of $n$ pairwise unlinked circles form the trivial link of $n$ components.

The ring group $R_n$ turns out to be closely related to several other groups that have been studied before in a variety of contexts under different names. This connection arises from one of our main technical results:

\begin{main}
The inclusion of $\RR_n$ into the space $\LL_n$ of all smooth trivial links of $n$ components in $\R^3$ is a homotopy equivalence.
\end{main}

Thus $R_n$ is isomorphic to the group $\pi_1\LL _n$ first studied in the 1962 thesis of Dahm [D], who identified it with a certain subgroup of the automorphism group of a free group on $n$ generators, subsequently called the symmetric automorphism group [Mc], [C]. A finite-index subgroup of this group is the `braid-permutation group' of [FRR]. Other references are [G], [BL], [R], [BMMM], [JMM], [BWC].

\begin{figure}[h!]
\begin{center}
\includegraphics[width=4.5in]{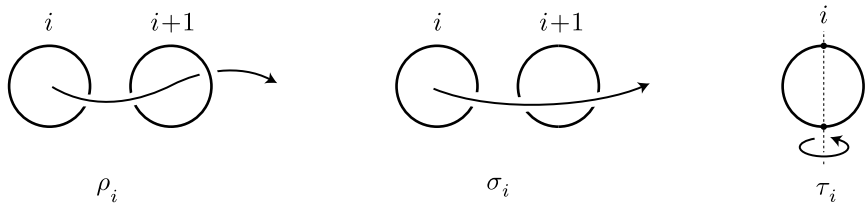}
\caption{ The generators $\rho_i$, $\sigma_i$, and $\tau_i$.}
\label{fig1}
\end{center}
\end{figure}

We will show that the group $R_n$ is generated by three families of elements $\rho_i$, $\sigma_i$, and $\tau_i$ defined as follows. If we place the $n$ rings in a standard position in the $yz$-plane with centers along the $y$-axis, then there are two generators $\rho_i$ and $\sigma_i$ that permute the $i$th and $(i+1)$st rings by passing the $i$th ring either through the $(i+1)$st ring or around it, respectively, as in Figure~\ref{fig1}.

\noindent
The generator $\tau_i$ reverses the orientation of the $i$th ring by rotating it 180 degrees around its vertical axis of symmetry. It is not hard to see that $\tau_i$ has order two in~$R_n$. We will show that the $\rho_i$'s generate a subgroup of $R_n$ isomorphic to the braid group $B_n$ and the $\sigma_i$'s generate a subgroup isomorphic to the symmetric group $\Sigma_n$.

\p{Parallel rings.} The space $\RR_n$ has a number of interesting subspaces.  The first of these we single out is the `untwisted ring space' $\UR_n$ consisting of all configurations of rings lying in planes parallel to a fixed plane, say the $yz$-plane. The loops of configurations giving the generators $\rho_i$ and $\sigma_i$ lie in this subspace. We will show that the untwisted ring group ${\ur}_n=\pi_1\,\UR_n$ is generated by the $\rho_i$'s and $\sigma_i$'s, and that the map ${\ur}_n\to R_n$ induced by the inclusion $\UR_n\incl\RR_n$ is injective, so ${\ur}_n$ can be identified with the subgroup of $R_n$ generated by the $\rho_i$'s and $\sigma_i$'s. We will also see that ${\ur}_n$ can be described as the fundamental group of the $2^n$\nobreakdash-sheeted covering space of $\RR_n$ consisting of configurations of oriented rings, so ${\ur}_n$ has index $2^n$ in~$R_n$. The $\tau_i$'s generate a complementary subgroup isomorphic to $\Z_2^n$, but neither this subgroup nor ${\ur}_n$ is normal in $R_n$.

Intermediate between $\UR_n$ and $\RR_n$ is the space $\VR_n$ of configurations of rings lying in vertical planes. Its fundamental group $\vr_n$ is also generated by the $\rho_i$'s, $\sigma_i$'s, and $\tau_i$'s, but the $\tau_i$'s have infinite order in $\vr_n$.

\p{Wickets.} Another interesting subspace of $\RR_n$ consists of configurations of rings, each of which is vertical and is cut into two equal halves by the $xy$-plane.  The upper halves of these rings can be thought of as wickets, as in the game of croquet, in upper half-space $\R^3_+$, and this subspace of $\RR_n$ can be identified with the space $\W_n$ of all configurations of $n$ disjoint wickets in $\R^3_+$. The condition of being pairwise unlinked is automatically satisfied for vertical rings that are bisected by the $xy$-plane.  In analogy to Theorem~1, one can compare $\W_n$ with the space $\A_n$ of configurations of $n$ disjoint smooth unknotted and unlinked arcs in $\R^3_+$ with endpoints on $\partial \R^3_+ = \R^2$. Here `unknotted and unlinked' means `isotopic to the standard configuration of $n$ disjoint wickets'.

\begin{main}
\label{main} 
The inclusion $\W_n \incl \A_n$ is a homotopy equivalence.
\end{main}

In fact, we will prove a common generalization of this result and Theorem 1 that involves configurations of both rings and wickets.

We call the group $\pi_1\W_n$ the {\it wicket group\/} $W_n$.  It too is generated by the $\rho_i$'s, $\sigma_i$'s, and $\tau_i$'s. The $\rho_i$'s again generate a subgroup isomorphic to $B_n$, but the $\sigma_i$'s now generate a subgroup that is isomorphic to $B_n$ rather than $\Sigma_n$.  The $\tau_i$'s have infinite order just as they do in $\vr_n$.  There is also an untwisted wicket group $UW_n=\pi_1\,\UW_n$ where $\UW_n=\W_n\cap \,\UR_n$.  We show that $UW_n$ is generated by the $\rho_i$'s and $\sigma_i$'s, and that the map $UW_n\to W_n$ induced by inclusion is injective, so $UW_n$ can be identified with the subgroup of $W_n$ generated by the $\rho_i$'s and $\sigma_i$'s. 

When defining ${\ur}_n$, $\vr_n$, $W_n$, and $UW_n$ as fundamental groups we did not mention basepoints, and this is justified by the fact that $\UR_n$, $\VR_n$, $\W_n$, and $\UW_n$ are all connected, by the same argument that shows that $\RR_n$ is connected. 

\smallskip
Summarizing, we have the following commutative diagram relating the various ring and wicket groups:

\begin{center}
\includegraphics[width=2in]{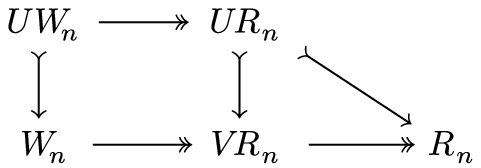}
\end{center}

\noindent
The two vertical maps are injective and correspond to adjoining the generators $\tau_i$. We will show that the two maps from the first column to the second column are quotient maps obtained by adding the relations $\sigma_i^2 = 1$, and the lower right horizontal map is the quotient map adding the relations $\tau_i^2=1$.

\p{Presentations} In Section~\ref{presentations}, we will derive finite presentations for all five of the groups in the diagram above, with the $\rho_i$'s, $\sigma_i$'s, and $\tau_i$'s as generators. The relations that hold for all five groups are the usual braid relations among the $\rho_i$'s and $\sigma_i$'s separately, together with certain braid-like relations combining $\rho_i$'s and $\sigma_i$'s, and for the groups in the second row there are relations describing how the $\tau_i$'s interact with the other generators.  For the three ring groups there are also the relations $\sigma_i^2=1$, and in $R_n$ the relations $\tau_i^2=1$ are added.

For ${\ur}_n$ the presentation was known previously [FRR], [BWC] using one of the more classical definitions of this group.  A presentation for $W_n$ was derived in [T1], using its interpretation as $\pi_1\A_n$, after generators had been found earlier in [H3].

\p{Asphericity} The space of configurations of $n$ distinct points in $\R^2$ is aspherical, with trivial higher homotopy groups, but this is no longer true for the ring spaces $\RR_n$, $\UR_n$, and $\VR_n$.  This is because the groups $R_n$, ${\ur}_n$, and $\vr_n$ contain torsion, the subgroup $\Sigma_n$ generated by the $\sigma_i$'s, so any $K(\pi,1)$ complex for these groups has to be infinite dimensional, but the spaces $\RR_n$, $\UR_n$, and $\VR_n$ are smooth finite-dimensional manifolds, hence finite-dimensional CW complexes (as are $\W_n$ and $\UW_n$). The situation is better for the wicket spaces:

\begin{main}
\label{two}
The spaces $ \W_n $ and $\UW_n$ are aspherical.
\end{main} 

In particular, this implies that $W_n$ and $UW_n$ are torsionfree. The proof of this theorem in Section~5 is more difficult than the proof of the corresponding result for configurations of points in $\R^2$, as it uses Theorem~2 as well as some results from $3$-manifold theory.

\p{Wicket groups as subgroups of braid groups.}  There is a natural homomorphism $W_n \to B_{2n}$ induced by the map which associates to each configuration of $n$ wickets the $2n$ endpoints of these wickets, a configuration of $2n$ points in $\R^2$. For example, the generators $\rho_i$ and $\sigma_i$ give rise to the two braids shown in Figure~2.

\begin{figure}[h!]
\begin{center}
\includegraphics[width=3.2in]{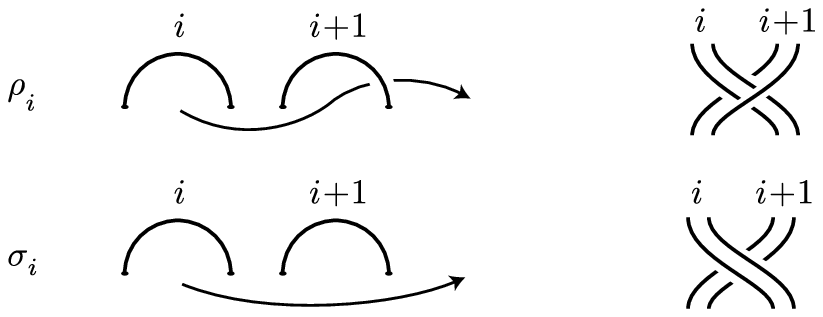}
\caption{ }
\label{fig2 }
\end{center}
\end{figure}

\noindent
It is a classical fact, whose proof we recall in Section 5, that this homomorphism $W_n \to B_{2n}$ is injective. This gives an alternative way of looking at $W_n$ as `braids' of $n$ ribbons, where certain intersections of ribbons are permitted, intersections that are known in knot theory as ribbon intersections.

\p{Pure versions.}  Just as the braid group $B_n$ has a pure braid subgroup $PB_n$, so do the five groups in the earlier commutative diagram have `pure' subgroups, the kernels of natural homomorphisms to $\Sigma_n$ measuring how loops of configurations permute the rings or wickets.  As in the braid case, these pure ring and wicket groups are the fundamental groups of the corresponding configuration spaces of ordered $n$-tuples of rings or wickets.  The full ring group $R_n$ is the semidirect product of the pure untwisted ring group ${\pur}_n$ and the signed permutation group $\Sigma^\pm_n$.  For the wicket group $W_n$ there is a weaker result, a nonsplit short exact sequence $0\to PUW_n\to W_n\to \Sigma^{\smZ}_n\to 0$ where $\Sigma^{\smZ}_n$ is the semidirect product of $\Sigma_n$ and $\Z^n$.

Our simple geometric method for finding presentations of the five `impure' ring and wicket groups also gives presentations for the pure versions of the ring groups $R_n$, ${\ur}_n$, and $\vr_n$, but not for the pure wicket groups. In the case of the pure untwisted ring group $\pur_n$ the presentation was originally found in [Mc]. It has generators $\alpha_{ij}$ in which all rings except the $i$th ring are stationary and the $i$th ring is pulled through the $j$th ring and back to its initial position without passing through any other rings, for each pair $i\ne j$.

\p{Rings of unequal sizes.} The subgroup of ${\pur}_n$ generated by the $\alpha_{ij}$ with $i<j$ has been studied in [CPVW].  We show that this `upper triangular pure untwisted ring group' is the fundamental group of the subspace $\UR^<_n$ of $\UR_n$ consisting of configurations of rings of unequal diameters. The sizes of the rings then provide a canonical ordering of the rings, hence loops in this space give elements of ${\pur}_n$, and we show the resulting homomorphism $\pi_1\,\UR^<_n \to {\pur}_n$ is injective with image the subgroup generated by the $\alpha_{ij}$ with $i<j$.

\p{Passing to the 3-sphere.}  In Section~\ref{sphere}, we also obtain similar results for configurations of circles in $S^3$ and wickets in a ball.  In the latter case wickets can be viewed as geodesics in hyperbolic $3$-space, and the configuration space of disjoint wickets is a subspace of the space of configurations of disjoint geodesics, a dense subspace having the same homotopy type as the larger space.  (A pair of disjoint geodesics can have an endpoint in common, so the two spaces are not identical.)

\p{Complexity of configurations.}
A key step in proving these results is a process for simplifying configurations of rings in $\RR_n$.  General configurations in $\RR_n$ can be quite complicated, with all the rings tightly packed together.  This happens already in the subspace $\W_n$ where the unlinking condition is automatic. One can take an arbitrary finite set of wickets, possibly intersecting in very complicated ways and tightly packed together, and then with a small random perturbation remove all the intersections to produce a configuration in $\W_n$.  The goal of the simplification process is to produce configurations in which each circle is surrounded by a region in which it is much larger than all other circles that intersect the region.  This region, or `microcosm', is by definition a closed ball of double the radius of the circle, and with the same center.  We define the complexity of a configuration of circles $C_1,\cdots,C_n$ of radii $r_1,\cdots,r_n$ to be the maximum of the ratios $r_i/r_j \le 1$ for the pairs of circles $C_i,C_j$ whose microcosms intersect. If none of the microcosms intersect, the complexity is defined to be $0$.  If we let $\RR^c_n$ be the subspace of $\RR_n$ consisting of configurations of complexity less than $c$, then the simplification process will show that the inclusion of $\RR^c_n$ into $\RR_n$ is a homotopy equivalence for any $c>0$.  

Configurations of small complexity can be thought of not only on the small scale of microcosms, but also in large-scale astronomical terms.  When the microcosms of two circles intersect, one can think of the smaller circle as a ring-shaped planet with the larger circle as its ring-shaped sun.  There can be several such planets in each solar system, each planet can have its own system of moons, the moons can have their own `moonlets,' and so on.  The solar systems can form galaxies, etc. 

The process of deforming $\RR_n$ into $\RR^c_n$ is an elaboration on the argument for showing $\RR_n$ is path-connected by shrinking all circles simultaneously in a canonical way.  If one starts with a configuration which is in general position in the sense that no circle has its center on the disk bounded by another circle, then this shrinking process produces a configuration of circles lying in disjoint balls.  This suffices to show $\RR_n$ is path-connected, but to capture its full homotopy type one cannot restrict attention to configurations that are in general position.  We deal with general configurations by combining shrinking with a pushing process that is realized by extending shrinkings of circles to ambient isotopies.  This is explained in detail in Section~\ref{complexity} of the paper.

\p{Configurations of spheres and disks.} The proof of Theorems 1 and 2, that the inclusions $\RR_n\incl\LL _n$ and $\W_n\incl\A_n$ are homotopy equivalences, uses the complexity reduction result described above, and it also involves a shift in focus from codimension two objects to codimension one objects, embedded spheres and disks, which are generally more tractable.  In Section~\ref{rigidifying} we use a parametrized disjunction technique to create the necessary configurations of spheres and disks, then we use the analogs of Theorems 1 and 2 for spheres and disks to improve configurations of smooth spheres and disks to round spheres and disks.  This relies ultimately on the proof of the Smale Conjecture in [H1], as does the final step of turning smooth circles and arcs into round circles and arcs. The spheres and disks are introduced to reduce the problem from configurations of many circles and arcs to configurations of at most one circle or arc in each complementary region of a configuration of spheres and disks.

\p{Dimension.} The paper concludes with a brief discussion in Section~\ref{dimension} of some elementary things that can be said about the homological dimension of the ring and wicket groups.

\section{Reducing Complexity.}
\label{complexity}

One way to define the topology on $\RR_n$ is in terms of its covering space consisting of ordered $n$-tuples of disjoint oriented circles in $\R^3$.  This covering space can be identified with an open subset in $\R^{6n}$ by assigning to each circle its centerpoint together with a vector orthogonal to the plane of the circle, of length equal to the radius of the circle and oriented according to the orientation of the circle via some rule like the right-hand rule. Ignoring ordering and orientations of circles amounts to factoring out the free action of the signed permutation group on this space.  Thus we see that $\RR_n$ has a finite-sheeted covering space which is an open set in $\R^{6n}$, and so $\RR_n$ itself is an open manifold of dimension $6n$. By similar reasoning one sees that the subspaces $\VR_n$, $\UR_n$, $\W_n$, and $\UW_n$ of $\RR_n$ are submanifolds of dimensions $5n$, $4n$, $4n$, and $3n$, respectively.

Let us recall the definition of complexity from the introduction.  If $C$ is a configuration in $\RR_n$ consisting of disjoint circles $C_1,\cdots,C_n$, let $B_i$ be the closed ball containing $C_i$ having the same center and double the radius.  (There is nothing special about the factor of $2$ here, and any other number greater than $1$ could be used instead.)  Then the {\it complexity\/} of the configuration $C$ is the maximum of the ratios $r_i/r_j\le 1$ of the radii of the pairs of circles $C_i, C_j$ in $C$ such that $B_i \cap B_j $ is nonempty, with the complexity defined to be $0$ if no $B_i$'s intersect.  We remark that complexity, as a function $\RR_n\to [0,1]$, is upper semicontinuous, meaning that small perturbations of a configuration $C$ cannot produce large increases in the complexity.  They can however produce large decreases if two circles $C_i,C_j$ whose balls $B_i,B_j$ intersect in a single point are perturbed so that $B_i$ and $B_j$ become disjoint.  

Define $\RR^c_n$ to be the subspace of $\RR_n$ consisting of configurations of complexity less than $c$.  This is an open subset of $\RR_n$.

\begin{thm} \label{1.1} The inclusion $ \RR^c_n \incl \RR_n $ is a homotopy equivalence for each $c>0$. The same is true for the subspaces $\UR^c_n\incl \UR_n$, $\VR^c_n\incl \VR_n$, $\W^c_n\incl \W_n$, and $\UW^c_n\incl \UW_n$.
\end{thm}

\vspace{-9pt}

\proof First we describe the argument from [FS], Lemma~3.2, for showing that $\RR_n$ is connected. Each configuration of disjoint circles in $ \R^3 $ bounds a unique configuration of hemispheres in $ \R^4_+ $ orthogonal to $ \R^3 $. The claim is that these hemispheres will be disjoint when each pair of circles is unlinked. To see this, think of $\R^4_+$ as the upper halfspace model of hyperbolic $4$-space, with the hemispheres as hyperbolic planes.  If two such planes intersect, they do so either in a single point or in a hyperbolic line, but the latter possibility is ruled out by the disjointness of the original collection of circles. Switching to the ball model of hyperbolic space, the point of intersection of two hyperbolic planes can be moved to the center of the ball, so the planes become Euclidean planes through the origin, and any transverse pair of such planes can be deformed through transverse planes to be orthogonal, when it is obvious that their boundary circles are linked in the boundary sphere~$S^3$. Thus unlinked circles in $\R^3$ bound disjoint hemispheres in $\R^4_+$.

For a configuration of circles in $\R^3$ bounding disjoint hemispheres in $\R^4_+$, consider what happens when one intersects the configuration of hemispheres with the hyperplanes $ \R^3_u = \R^3 \times \{ u \} $ for $ u \ge 0 $. As $ u $ increases, each circle shrinks to its centerpoint and disappears. Let us call this the {\em canonical shrinking\/} of the configuration. 

A given configuration of circles can be perturbed so that no centerpoint of one circle lies on the disk bounded by another circle. Then if we perform the canonical shrinking of the configuration, we can stop the shrinking of each circle just before it shrinks to a point and keep it at a small size so that no other shrinking circles will bump into it. In this way the given circle configuration can be shrunk until the disks bounded by the circles are all disjoint. This says that the configuration of circles forms the trivial link, and it makes clear that the space $ \RR_n $ is path-connected. 

When dealing with a $k$-parameter family of circle configurations, however, one cannot avoid configurations where one circle has center lying in the disk bounded by another circle.  If the latter circle is larger than the first, the two circles would then collide if we stop the shrinking of the smaller circle just before it disappears.  Our strategy to avoid such collisions will still be to stop the canonical shrinking of each circle just before it disappears, and thereafter shrink it at a slower rate so that it does not disappear, but we also allow it to be pushed by `air cushions' surrounding larger circles as they shrink, so that the smaller circle never intersects the larger circles. 

The pushing will be achieved by an inductive process that relies on extending isotopies of circles to ambient isotopies of $\R^3$, so let us recall the standard procedure in differential topology for extending isotopies of submanifolds to ambient isotopies.  An isotopy of a submanifold $N$ of a manifold $N$ is a level-preserving embedding $F\cln N\times I \incl M\times I$.  This has a tangent vector field given by the velocity vectors of the paths $t\mapsto F(x,t)$.  The second coordinate of this vector field is equal to $1$, and we extend it to a vector field on $M\times I$ with the same property by damping off the first coordinate to $0$ as one moves away from $F(N\times I)$ in a small tubular neighborhood of $F(N\times I)$.  Then the flow lines of this extended vector field define the extended isotopy.  This also works with $I$ replaced by $[0,\infty)$ as will be the case in our situation.  The manifold $M$ will be $\R^3$, and we can choose the tubular neighborhood of the submanifold $F(N\times [0,\infty))$ to be an $\epsilon(t)$-neighborhood of $F(N\times \{t\})$ in each level $\R^3\times\{t\}$.  

For a configuration $C$ in $\RR_n$ consisting of circles $C_1,\cdots,C_n$, let $C^1$ be the union of the largest circles in $C$, $C^2$ the union of the next-largest circles, and so on.  Let $u$ be the time parameter in the canonical shrinking of $C$, and let $u=u_i$ be the time when the circles of $C^i$ shrink to their centerpoints, so $u_1 > u_2> \cdots$ . Note that all the circles in $C^i$ have distinct centerpoints since two circles with the same center and radius must intersect.  The canonical shrinking defines an isotopy $\Phi_u(C^i)$ for $u<u_i$.  Our aim is to truncate this at a value $u=u'_i$ slightly less than $u_i$, then extend the truncated isotopy to values of $u$ greater than $u'_i$.  The new extended isotopy $\Phi_u(C^i)$ will move each circle $C_j$ of $C^i$ through circles parallel to itself, so $\Phi_u(C_j)$ will be determined by specifying the centerpoint $c_j(u)$ and the radius $r_j(u)$ of $\Phi_u(C_j)$.  The center $c_j(u)$ is the centerpoint of $C_j$ for $u\le u'_i$ since this point does not move during the canonical shrinking, and we will in fact have $c_j(u)$ equal to this same point for $u\le u_i$, not just $u\le u'_i$.  For the function $r_j(u)$, the canonical shrinking specifies this for $u\le u'_i$, and we will choose it to be a positive decreasing function of $u$ for $u>u'_i$.

The extended isotopy $\Phi_u(C^i)$ will be constructed by induction on $i$.  For $i=1$ and $C_j$ a circle of $C^1$ we let $c_j(u)$ be constant for all $u$, and we let $r_j(u)$ be any decreasing function $r^1(u)$ of $u$ for $u > u'_1$ where $u'_1$ is chosen close enough to $u_1$ so that the microcosms of all the circles of $\Phi_{u'_1}(C^1)$ are disjoint.  Such a $u'_1$ exists since the centerpoints of the circles of $C^1$ are distinct.  The microcosms of the circles of $\Phi_u(C^1)$ will then remain disjoint for all $u>u'_1$. To finish the first step of the induction we extend the isotopy $\Phi_u(C^1)$ to an ambient isotopy $\Phi^1_u\cln \R^3\to \R^3$ by the general procedure described earlier, with $\Phi^1_0$ the identity.

For a circle $C_j$ of $C^2$ with centerpoint $c_j$ we let $c_j(u)$ be constant for $u\le u_2$ and then we let it move via the isotopy $\Phi^1_u$.  In formulas this means $c_j(u) = \Phi^1_u(\Phi^1_{u_2})^{-1}(c_j)$.  This will in fact be constant for $u$ slightly greater than $u_2$ as well as for $u\le u_2$.  Since $\Phi^1_u$ is an ambient isotopy, $c_j(u)$ will be disjoint from $\Phi_u(C^1)$ and from $c_k(u)$ for other circles $C_k$ of $C^2$ for all $u$.  This implies that if we choose $u_2'$ close enough to $u_2$ and the function $r^2(u)$ giving the radius of the circles of $\Phi_u(C^2)$ small enough, then these circles will be disjoint from $\Phi_u(C^1)$ for all $u$ and will have disjoint microcosms for $u>u'_2$.  We can also make $r^2(u)$ small enough so that the ratio $r^2(u)/r^1(u)$ goes to $0$ with increasing $u$.  The second step of the induction is completed by extending the isotopies $\Phi_u(C^1)$ and $\Phi_u(C^2)$ to an ambient isotopy $\Phi^2_u$ starting with $\Phi^2_0$ the identity.  

Subsequent induction steps are similar.  For example, at the next stage, for a circle $C_j$ of $C^3$ with centerpoint $c_j$ we let $c_j(u)$ move according to the isotopy $ \Phi^2_u$, and we choose $u'_3$ close enough to $u_3$ and $r^3(u)$ small enough so that the resulting circles of $\Phi_u(C^3)$ are disjoint from $\Phi_u(C^1)$ and $\Phi_u(C^2)$ for all $u$ and the microcosms of the circles of $\Phi_u(C^3)$ are disjoint for $u>u'_3$.  Also we make $r^3(u)$ small enough so that the ratio $r^3(u)/r^2(u)$ goes to $0$ with increasing $u$. We can also assume that $r^3(u)/r^2(u)<r^2(u)/r^1(u)$, and inductively that $r^{i+1}(u)/r^i(u) < r^i(u)/r^{i-1}(u)$ for all $i$.

When the induction process is finished we have a path $\Phi_u(C)$ in $\RR_n$, defined for each $C\in \RR_n$.  It is clear that the complexity of $\Phi_u(C)$ approaches $0$ as $u$ goes to $\infty$ since the circles of $\Phi_u(C^i)$ have disjoint microcosms for large $u$ and the ratios $r^{i+1}(u)/r^i(u)$ approach $0$.  We claim that the complexity of $\Phi_u(C)$ decreases monotonically (in the weak sense) as $u$ increases.  Consider two circles of $C$, say $C_1$ and $C_2$.  If they are in the same $C^i$, they have the same radius throughout the isotopy $\Phi_u$, and their centers are stationary until $u=u_i$, after which their microcosms remain disjoint, so their contribution to the complexity decreases monotonically, being either $0$ for all $u$ or $1$ for a while and then $0$.  If $C_1$ and $C_2$ belong to different $C^i$'s, with $C_1$ in $C^{i_1}$ and $C_2$ in $C^{i_2}$ for $i_1>i_2$, the ratio of their radii approaches $0$ monotonically, so the only way they could contribute to a non-monotonic complexity would be for their microcosms to bump into each other at a certain time $u$ after having been disjoint shortly before this time.  For this to happen, both $\Phi_u(C_1)$ and $\Phi_u(C_2)$ would have to be within the microcosm of some larger circle $\Phi_u(C_3)$ in $C^{i_3}$ for some $i_3<i_2$.  In this case the pair $\Phi_u(C_2)$, $\Phi_u(C_3)$ would be contributing a larger number to the complexity than the pair $\Phi_u(C_1)$, $\Phi_u(C_2)$, so the collision between the microcosms of the latter pair would not be causing an increase in the overall complexity.

To show that the inclusion $\RR^c_n\incl \RR_n$ is a homotopy equivalence for $c>0$  it suffices to show that the relative homotopy groups $\pi_k(\RR_n,\RR^c_n)$ are zero for all $k$, since both spaces are smooth manifolds and hence CW complexes.  Thus it suffices to deform a given a map $(D^k,\bdy D^k)\to (\RR_n,\RR^c_n)$, $t\mapsto C_t$, through such maps to a map with image in $\RR^c_n$. This would follow if we could add a parameter $t\in D^k$ to our previous construction of the deformation $\Phi_u$.  However, there is a problem with doing this directly because the relative sizes of the circles in a family of configurations $C_t \in \RR_n$ can change with varying $t$, so the sequence of induction steps in the construction of the desired deformation $\Phi_{tu}$ could change with $t$.  What we will do instead is concatenate initial segments of deformations $\Phi_{tu}$ over different regions in $D^k$ to produce a new family of deformations $\Psi_{tu}$.

As a preliminary step, note that choosing an ordering of the circles of the configuration $C_t$ for one value of $t$ gives an ordering for all $t$ since the parameter domain $D^k$ is simply-connected.  Thus we can label the circles as $C^t_1,\cdots,C^t_n$.  The radius of $C^t_i$ varies continuously with $t$, and we can approximate these radius functions arbitrarily closely by piecewise linear functions of $t$, close enough so that they correspond to a deformation of the family $C_t$, staying in the open set $\RR^c_n$ over $\bdy D^k$.  Thus we may assume the radius functions are piecewise linear.  This means we can triangulate $D^k$ so that the radius functions are linear on simplices.  After a subdivision of this triangulation, we can assume that on the interior of each simplex the ordering of the circles $C^t_i$ according to size is constant, and as one passes to faces of a simplex all that happens to this ordering is that some inequalities among sizes become equalities.  

We will construct the final deformations $\Psi_{tu}$ by a second induction, where the inductive step is to extend $\Psi_{tu}$ from a neighborhood of the $p$-skeleton of the triangulation of $D^k$ to a neighborhood of the $(p+1)$-skeleton.  More specifically, we will construct continuous functions $\psi_0\le \psi_1\le\cdots\le\psi_k$ from $D^k$ to $[0,\infty)$ such that the inductive step consists of extending $\Psi_{tu}$ from being defined for $0\le u\le\psi_p(t)$ to being defined for $0\le u\le\psi_{p+1}(t)$. The functions $\psi_p$ will satisfy:

\begin{list}{}{\setlength{\leftmargin}{22pt}\setlength{\labelwidth}{16pt}\setlength{\labelsep}{5pt}\setlength{\itemsep}{3pt}}

\item[(a)]  $\psi_p = 0$ outside some neighborhood $N_p$ of the $p$-skeleton.
%\smallskip
\item[(b)] $\Psi_{tu}(C_t)$ lies in $\RR^c_n$ for $u=\psi_p(t)$ when $t$ lies in a smaller neighborhood $N'_p$ of the $p$-skeleton.
%\smallskip
\item[(c)] $\psi_p =\psi_{p+1}=\cdots =\psi_k$ in $N'_p$. 

\end{list}
%\smallskip\noindent
The ordering of the circles of $C_t$ according to size will be preserved during the deformation $\Psi_{tu}$.

For the induction step of extending $\Psi_{tu}$ over a $p$-simplex $\sigma$, let $\sigma'$ be a slightly smaller copy of $\sigma$ lying in the interior of $\sigma$ and with boundary in the interior of $N'_{p-1}$.  As $t$ varies over $\sigma'$ the size ordering of the circles of $C_t$ is constant. For each $t$ in $\sigma'$ we apply the earlier inductive procedure to construct a deformation $\Phi_{tu}$, starting with the family $\Psi_{tu}(C_t)$ for $u=\psi_{p-1}(t)$.  This can be done continuously in $t\in\sigma'$ since the various choices in the construction can be made to vary continuously with $t$.  These choices are: the numbers $u'_i(t)<u_i(t)$, the radius functions $r^i(t,u)$, and the isotopy extensions $\Phi^i_{tu}$. The construction of $\Phi_{tu}$ works in fact in a neighborhood of $\sigma'$ in $D^k$ by extending the functions $u'_i(t)$ and $r^i(t,u)$ and the isotopy extensions $\Phi^i_{tu}$ to nearby $t$ values.  As $t$ moves off $\sigma'$ the size ordering in $C_t$ may vary, as some size equalities become inequalities, but we still use the same decomposition of $C_t$ into the subsets $C^i_t$, and we choose the functions $r^i(t,u)$ so that for each $t$ in the neighborhood, this size ordering is preserved throughout the deformation $\Phi_{tu}$. To finish the induction step we choose $\psi_p$ by requiring $\psi_p - \psi_{p-1}$ to have support in a neighborhood of $\sigma'$ and to have large enough values in a smaller neighborhood of $\sigma'$ so that $\Phi_{tu}(C_t)$ lies in $\RR^c_n$ for $t$ in this smaller neighborhood and $u\ge \psi_p(t)$. Then we extend the previously defined $\Psi_{tu}(C_t)$ for $u\in [0,\psi_{p-1}(t)]$ by defining it to be equal to $\Phi_{tu}(C_t)$ for $u\in [\psi_{p-1}(t),\psi_p(t)]$.

This finishes the proof for the inclusion $\RR^c_n\incl \RR_n$.  Since the deformations $\Phi_{tu}$ take circles to parallel circles, the proof also applies for the inclusions $\UR^c_n\incl \UR_n$ and $\VR^c_n\incl \VR_n$.  For the inclusions $\W^c_n\incl \W_n$ and $\UW^c_n\incl \UW_n$, observe that in the case of configurations of wickets, the extended isotopies $\Phi^i_{tu}$ take the $xy$-plane to itself so they take wickets to wickets.  \qed

\smallskip

\p{Remarks  on the proof of Theorem~\ref{1.1}.} We can strengthen the proof slightly to give a deformation of the given family $C_t$ to a family which not only has small complexity but has the additional property that the microcosm around each circle is disjoint from all larger circles.  This can be achieved by choosing the radius function $r^i(u)$ sufficiently small at each stage of the construction of the deformations $\Phi_u$.  In the later part of the proof when $\Psi_{tu}$ is constructed from truncated deformations $\Phi_{tu}$, initial segments of canonical shrinkings are also inserted, and these preserve the additional property since smaller circles shrink faster than larger circles during the canonical shrinking.

The proof also works for the configuration space $\WR_{m,n}$ consisting of configurations of $m$ wickets and $n$ rings in $\R^3_+$, all the wickets and rings being disjoint and pairwise unlinked, and with the rings disjoint from the $xy$-plane.  Thus $\WR_{m,0}$ is $\W_m$, and it is easy to see that $\WR_{0,n}$ and $\RR_n$ are homeomorphic, although they are not identical since one consists of configurations in $\R^3_+$ and the other of configurations in $\R^3$. Namely, both contain the space of configurations of rings for which the minimum $z$-value of all the rings is $1$, and $\WR_{0,n}$ is the product of this subspace with $(0,\infty)$ while $\RR_n$ is the product of this subspace with $\R$. In each case the projection onto the first factor is given by vertically translating configurations to make their minimum $z$-value $1$, and projection onto the second factor is by taking the minimum $z$-value of a configuration.

\p{A further enhancement.} A slight variation on the technique used to prove the theorem will be used to prove the following result:

\begin{prop} 
\label{1.2}
The natural maps ${\ur}_n\to R_n$, ${\ur}_n\to \vr_n$, and $ UW_n \to W_n $ induced by the inclusions $\UR_n\incl\RR_n$, $\UR_n\incl \VR_n$, and $ \UW_n \incl \W_n $ are injective.
\end{prop}

\vspace{-9pt}

\proof Consider first the case of ${\ur}_n\to R_n$.  Let $\PUR_n$ and $\PR_n$ be the ``pure" versions of $\UR_n$ and $\RR_n$, the covering spaces of $\UR_n$ and $\RR_n$ obtained by ordering the rings, so that $\UR_n$ and $\RR_n$ are the quotients of $\PUR_n$ and $\PR_n$ with the action of the symmetric group $\Sigma_n$ factored out. It will suffice to show injectivity of the map $\pi_1\PUR_n\to\pi_1\PR_n$ induced by the inclusion $\PUR_n\incl \PR_n$.

By associating to each ring in $\R^3$ the line through the origin orthogonal to the plane containing the ring we obtain a map $\PR_n\to (\RP^2)^n$ whose fibers over points in the diagonal of $(\RP^2)^n$ are copies of $\PUR_n$. Let us suppose for the moment that this map is a fibration.  It has a section, obtained by choosing a standard configuration of rings lying in disjoint balls and taking all possible rotations of these rings about their centers.  The existence of the section would then imply that the long exact sequence of homotopy groups breaks up into split short exact sequences, so in particular there would be a short exact sequence
$$
0\to \pi_1\,\PUR_n \to\pi_1\PR_n \to \pi_1(\RP^2)^n\to 0
$$
which would give the desired injectivity.

We will make this into a valid argument by showing the weaker result that the projection $\PR_n\to (\RP^2)^n$ is a quasifibration.  Recall that a map $p\cln E\to B$ is a quasifibration if $p_*\cln\pi_i(E,p^{-1}(b),e)\to\pi_i(B,b)$ is an isomorphism for each $b\in B$, $e\in p^{-1}(b)$, and $i\ge 0$.  Thus a quasifibration has a long exact sequence of homotopy groups just like for a fibration.  The standard argument for showing that a map $p\cln E\to B$ with the homotopy lifting property for maps of disks $D^k$, $k\ge 0$, has an associated long exact sequence of homotopy groups in fact proceeds by showing that the quasifibration property is satisfied.  This argument generalizes easily to a slightly weaker version of the homotopy lifting property, which asserts the existence of a lift, not of a given homotopy $D^k \times I\to B$, but of some reparametrization of this homotopy, obtained by composition with a map $D^k \times I \to D^k \times I$ of the form $(x,t)\mapsto (x,g_x(t))$ for a family of maps $g_x\cln (I,0,1)\to(I,0,1)$.  We will use this generalization below.

To show that the projection $\PR_n\to (\RP^2)^n$ is a quasifibration, the key observation is that we can enhance the construction of the deformations $\Phi_{tu}$ by not only shrinking the rings and moving their centers, but also rotating the rings according to any deformation of the planes that contain them, provided that we delay the start of these deformations to the time $u=u_1(t)$. At the inductive step when $\Phi^i_{tu}$ is constructed for the rings of $C^i_t$ for $u\ge u_i(t)$, these rings lie in microcosms that are disjoint from each other and from the larger rings for which $\Phi_{tu}$ has already been constructed, so they can be rotated arbitrarily about their centers, starting at time $u=u_1(t)$.  

With this elaboration on the construction of $\Phi_{tu}$ we construct the deformations $\Psi_{tu}$ as before.  First we deform a given map $D^k\to \PR_n$ to make the radii of the rings piecewise linear functions of the parameter $t\in D^k$.  Then we proceed by induction over the skeleta of the triangulation of $D^k$.  Prior to the induction step of extending over $p$-simplices, the deformation $\Psi_{tu}$ for $u\le\psi_{p-1}(t)$ will include some initial segment of a given deformation of the planes of the rings of $C_t$, reparametrized by the insertion of pauses. Then we construct $\Phi_{tu}$ as in the preceding paragraph, starting with $\Psi_{tu}(C_t)$ for $u = \psi_{p-1}(t)$. Thus the deformation of the planes containing the rings pauses for a time before continuing with the given deformation. At the end of the induction step we choose the function $\psi_p$ and truncate $\Phi_{tu}$, which can truncate the deformation of the planes containing the rings, so that they pause once more in the next stage of the induction.  It is no longer necessary to choose $\psi_p$ large enough to make $\Psi_{tu}(C_t)$ lie in $\RR^c_n$ for $u=\psi_p(t)$ if $t$ is near the $p$-skeleton.  Instead, we only need  it large enough to allow time to carry out the deformation of the planes of the rings.  

At the end of the induction process we have a deformation $\Psi_{tu}$ such that the planes of the rings vary by a reparametrization of the given deformation of these planes.  The parameter $u$ varies over an interval $[0,\psi_k(t)]$ but we can rescale to make this $[0,1]$. This finishes the proof that the projection $\PR_n\to (\RP^2)^n$ is a quasifibration, and hence the proof that ${\ur}_n\to R_n$ is injective. 

Since the injection ${\ur}_n\to R_n$ factors through $\vr_n$ it follows that ${\ur}_n\to \vr_n$ is also injective.  For $UW_n\to W_n$ we can use the same quasifibration argument as in the first case, the only difference being that $(\RP^2)^n$ is replaced by $(\RP^1)^n$, an $n$-dimensional torus.  \qed

\medskip

Another result stated in the introduction can be proved using the same method:

\begin{prop}
The natural map from $\UR_n$ to the covering space $\RR^+_n$ of $\RR_n$ consisting of configurations of oriented rings induces an isomorphism $\pi_1\,\UR_n \to \pi_1\RR^+_n$.
\end{prop}

\vspace{-9pt}
\proof  The arguments in the preceding proof work equally well with oriented rings, the only difference being that $\RP^2$ is replaced by $S^2$.  Since this is simply-connected, the previous short exact sequence of fundamental groups for the quasifibration reduces to an isomorphism $\pi_1\,\UR_n \to \pi_1\RR^+_n$.  \qed
\medskip

The short exact sequence 
$$
0\to {\pur}_n\to PR_n\to \Z_2^n \to 0
$$
constructed in the proof of Proposition \ref{1.2} has an obvious splitting, obtained by rotating the rings within disjoint balls.  This sequence embeds in a larger split short exact sequence
$$
0\to {\pur}_n\to R_n\to \Sigma^\pm_n\to 0
$$
where $\Sigma^\pm_n$ is the signed permutation group, the semidirect product of $\Sigma_n$ and $\Z^n_2$.  The homomorphism $R_n\to \Sigma^\pm_n$ assigns to each loop in $\RR_n$ the permutation of the rings that it effects, as well as the changes of orientations of the rings.  The sequence splits since $\Sigma^\pm_n$ is the fundamental group of the subspace of $\RR_n$ consisting of configurations of rings with disjoint microcosms.  This short exact sequence maps to another split exact sequence
$$
0\to PR_n\to R_n\to \Sigma_n\to 0
$$
which in turn contains the split exact sequence
$$
0\to {\pur}_n\to {\ur}_n\to \Sigma_n\to 0
$$
There are analogous sequences with $\vr_n$ in place of $R_n$ and with the $\Z_2$'s replaced by $\Z$'s and $\Sigma^\pm_n$ replaced by $\Sigma^{\smZ}_n$, the semidirect product of $\Sigma_n$ and $\Z^n$.  For wicket groups there are similar short exact sequences as well, but the only one that splits is the one not involving $\Sigma_n$, namely
$$
0\to PUW_n\to PW_n\to \Z^n \to 0
$$

\section{Presentations}\label{presentations}

In this section we use the results in the preceding section to obtain finite presentations of ring and wicket groups.  First an elementary result:
 
 \begin{prop}
 The elements $\sigma_i$ of $UW_n$ generate a subgroup isomorphic to the braid group $B_n$, and so also do the elements $\rho_i$.
 \end{prop}
 
 \vspace{-9pt}
 
 \proof  Let us take $ \UW_n $ to be the subspace of $ \W_n $ consisting of configurations of wickets lying in planes perpendicular to the $ x $-axis. Sending each wicket to its endpoint with larger $y$-coordinate defines a map $ \UW_n \to \C_n $ where $\C_n$ is the space of configurations of $n$ distinct points in $\R^2$, so $B_n=\pi_1\C_n$. The restriction of this map to the subspace $ \UW^{\sigma}_n $ of $ \UW_n $ consisting of configurations of wickets having disjoint projections to the $ xy $-plane is a homotopy equivalence. The maps $ \UW^{\sigma}_n \incl \UW_n \to \C_n $ induce homomorphisms $ B_n \to W_n \to B_n $ whose composition is the identity. The image of the first homomorphism is generated by the $\sigma_i$'s, so this subgroup of $ W_n $ is isomorphic to $ B_n $. 

The argument for $ \rho_i $'s is similar using the subspace $ \UW^{\rho}_n $ of $ \UW_n $ consisting of configurations of wickets, each of which is symmetric with respect to reflection across the $ xz $-plane. Wickets with this symmetry property are determined by their endpoints in the upper half of $ \R^2 $, so $ \UW^{\rho}_n $  can be identified with $ \C_n $ viewed as the space of configurations of $ n $ points in the upper half of $ \R^2 $.   \qed

\medskip

These arguments do not work with ${\ur}_n$ in place of $UW_n$, but the $\rho_i$'s still generate a copy of $B_n$ in ${\ur}_n$ as we will show in Proposition~4.2.  The $\sigma_i$'s, on the other hand, generate a copy of $\Sigma_n$ in ${\ur}_n$ since they have order $2$ and satisfy the braid relations, so the canonical map ${\ur}_n\to \Sigma_n$ has a section.

Now we determine a presentation for $ UW_n $ by a straightforward elaboration of the standard procedure for computing a presentation for $ B_n $ using general position arguments.

\begin{prop}
The group $UW_n$ has a presentation with generators the elements $\sigma_i$ and $\rho_i$ for $ i = 1, \cdots, n-1 $ and with the following relations:
%\smallskip
\begin{enumerate}
\item[]%{}
\hspace{-6pt}$[\rho_i,\rho_j] = [\sigma_i,\sigma_j] = [\rho_i,\sigma_j] = 1 $  if\hspace{3pt}  $ | i - j | > 1 $
\smallskip
\item[]%{}
\hspace{-6pt}$\rho_i\rho_{i+1}\rho_i = \rho_{i+1}\rho_i\rho_{i+1}$, \quad $\sigma_i\sigma_{i+1}\sigma_i = \sigma_{i+1}\sigma_i\sigma_{i+1}$
\smallskip
\item[]%{}
\hspace{-6pt}$\rho_i\sigma_{i+1}\sigma_i = \sigma_{i+1}\sigma_i\rho_{i+1}$, \quad 
   $\sigma_i\sigma_{i+1}\rho_i = \rho_{i+1}\sigma_i\sigma_{i+1}$, \quad 
$\sigma_i\rho_{i+1}\rho_i = \rho_{i+1}\rho_i\sigma_{i+1}$
\end{enumerate}
\end{prop}

\vspace{-7pt}

\proof We again take $ \UW_n $ to consist of configurations of wickets lying in planes perpendicular to the $ x $-axis. Let $ \UW^0_n $ be the open dense subspace of $ \UW_n $ consisting of configurations of wickets all lying in distinct planes. This subspace is homeomorphic to $ \R^{3n} $, so it is contractible. The complement of $ \UW^0_n $ decomposes into a disjoint union of connected manifold strata, determined by which subsets of wickets lie in the same planes and how these wickets are nested in these planes. Each stratum is homeomorphic to a Euclidean space of the appropriate dimension. The codimension one strata are formed by configurations with exactly two wickets lying in the same plane. These form a codimension one submanifold $ \UW^1_n $ of $ \UW_n $ defined locally by equating the $x$-coordinates of two wickets. The codimension two strata, forming a codimension two submanifold $ \UW^2_n $, consist of configurations where either two disjoint pairs of wickets lie in coinciding planes, or three wickets lie in a single plane. 

To find generators for $ UW_n $ consider a loop in $ \UW_n $. By general position this can be pushed off all strata of codimension $ 2 $ and greater until it lies in $ \UW^0_n \cup \UW^1_n $, and we may assume it is transverse to $ \UW^1_n $, crossing it finitely many times. Each such crossing corresponds to a generator $ \rho_i $ or $ \sigma_i $ or its inverse. Since the strata of $ \UW^1_n $ are contractible, they have trivial normal bundles and we can distinguish between the directions of crossing these strata. Since $ \UW^0_n $ is contractible, it follows that the given loop in $ \UW_n$ is homotopic to a product of $ \rho_i $'s and $ \sigma_i $'s and their inverses, so these elements generate $ UW_n $.

To find a complete set of relations among these generators, consider a homotopy in $ \UW_n $ between two loops of the type just considered. General position allows us to push this homotopy off strata of codimension greater than $ 2 $, and we can make it transverse to strata of $ \UW^2_n $ and $ \UW^1_n $. Let us examine what happens near points where the homotopy crosses $ \UW^2_n $. For strata of $ \UW^2_n $ where two disjoint pairs of wickets lie in coinciding planes we just have simple commuting relations: $ \rho_i $ and $ \sigma_i $ commute with $ \rho_j $ and $ \sigma_j $ if $ | i-j | > 1 $. More interesting are the relations arising from three wickets lying in the same plane. Here there are five cases according to how the projections of the wickets to the $ xy $-plane intersect. The three projections can be completely disjoint, completely nested, or some combination of disjoint and nested, as indicated in the first column of Figure~3, where for visual clarity we have perturbed the overlapping projections of the three wickets so that they appear to be disjoint.

\begin{figure}[h!] %[htp]
\begin{center}
\includegraphics[width=4in]{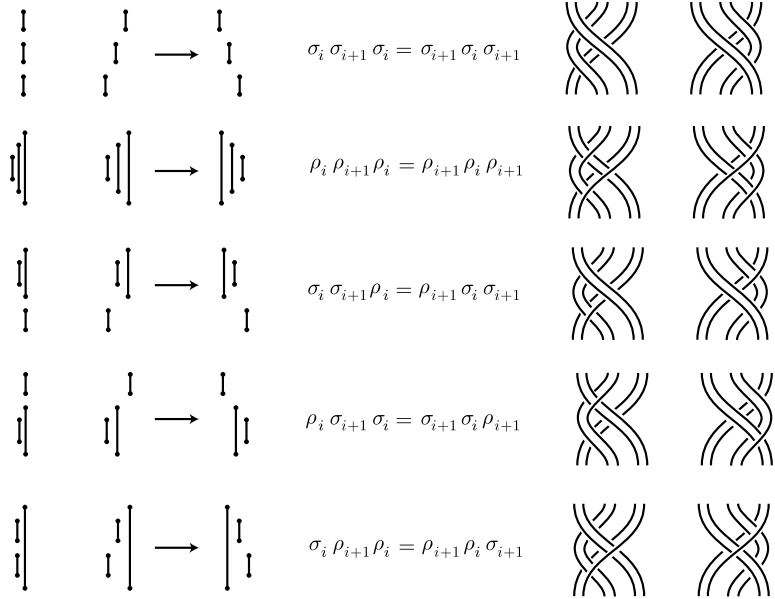}
\caption{Relations in $UW_n$. }
\label{rels}
\end{center}
\end{figure}

A small loop around the codimension $2$ stratum crosses codimension $1$ strata six times since the local picture is like the intersection of the three planes $ x = y $, $x = z $, and $ y = z $ in $ \R^3 $. One can view the resulting relation as an equation between two ways of going halfway around the codimension $2$ stratum. The starting and ending points of the two ways are shown in the second column of Figure~\ref{rels}. The relation itself is written in the next column, and the final column shows the braid picture of the relation, using the endpoint map $ UW_n \to B_{2n} $.   \qed

\medskip

\begin{prop}
A presentation for $ {\ur}_n $ is obtained from the presentation for $UW_n$ in the preceding proposition by adding the relations $\sigma_i^2 =1$.
\end{prop}

\vspace{-3pt}

Note that the relations $\rho_i\sigma_{i+1}\sigma_i = \sigma_{i+1}\sigma_i\rho_{i+1}$ and  $\sigma_i\sigma_{i+1}\rho_i = \rho_{i+1}\sigma_i\sigma_{i+1}$ in the presentation for $UW_n$ become equivalent if $\sigma_i$ and $\sigma_{i+1}$ have order $2$, so either relation can be omitted from the presentation for ${\ur}_n$.  The geometric explanation for this is that the third and fourth configurations in the preceding large diagram are obviously equivalent when we are dealing with rings rather than wickets.

\proof The argument is similar to that for $UW_n$.  We take $ \UR_n $ to consist of the configurations of rings lying in planes parallel to the $ xz $-plane. Strata here are defined just as for $ \UW_n $ according to the coincidences among these planes. The only essential difference is that now not all strata are contractible. A codimension one stratum where two planes coincide and the two rings in this plane are not nested has the homotopy type of a circle. Crossing this stratum corresponds to a generator $ \sigma_i $. The normal bundle of this stratum is nontrivial, which means that we cannot distinguish between $ \sigma_i $ and $ \sigma_i^{-1} $, or in other words, we have the relation $ \sigma_i^2 = 1 $. An alternative way to proceed would be to subdivide this stratum into two contractible codimension one strata separated by a codimension two stratum, the configurations where the centers of the two rings in this plane have the same projection to the $xy$-plane. A small loop around this codimension two stratum would give the relation $ \sigma_i^2 = 1 $. Using either approach we conclude that adding the relations $ \sigma_i^2 = 1 $ to the earlier presentation for $ UW_n $ gives a presentation for $ {\ur}_n $. \qed

\medskip

Next we turn to the pure untwisted ring group $\pur_n$. Recall the elements $\alpha_{ij}$ passing the $i$th ring through the $j$th ring and back to its initial position, for $i\ne j$.

\begin{prop}
The group $\pur_n$ has a presentation with generators the elements $\alpha_{ij}$ for $1\le i,j\le n$, $i\ne j$, and relations
$$
\alpha_{ij}\alpha_{k\ell}=\alpha_{k\ell}\alpha_{ij}\ \qquad
\alpha_{ik}\alpha_{jk}=\alpha_{jk}\alpha_{ik}\ \qquad
\alpha_{ij}\alpha_{ik}\alpha_{jk}=\alpha_{jk}\alpha_{ik}\alpha_{ij}
$$
where distinct symbols for subscripts denote subscripts that are distinct numbers.
\end{prop}

\vspace{-3pt}
Using the second relation, the third relation can be restated as saying that $\alpha_{jk}$ commutes with $\alpha_{ij}\alpha_{ik}$.

\vspace{-4pt}
\proof  The group ${\pur}_n$ is the fundamental group of the covering space $\PUR_n$ of $\UR_n$ in which the rings are numbered.  Let $\PUR^0_n$ be the subspace of $\PUR_n$ consisting of configurations in which no circles are nested within the planes that contain them.  We claim that $\PUR^0_n$ is simply-connected.  To see this, consider the projection of $\PUR^0_n$ to the space of ordered $n$-tuples of distinct points in $\R^3$ sending a configuration of circles to the configuration of its centerpoints. This projection has a section, sending a configuration of points to the configuration of circles of radius equal to one-quarter of the minimum distance between the points.  Further, $\PUR^0_n$ deformation retracts to the image of this section by first shrinking the circles whose radius is too large, then expanding the circles whose radius is too small.  Since the space of point configurations is simply-connected (by a standard induction argument involving fibrations obtained by forgetting one of the points), it follows that $\pi_1\,\PUR^0_n = 0$.  

Let $\PUR^1_n$ be obtained from $\PUR^0_n$ by adjoining the codimension-one strata, the configurations having exactly one circle nested inside another. The map $\pi_1\,\PUR^1_n\to \pi_1\,\PUR_n$ is surjective, so we see that ${\pur}_n$ is generated by the elements $\alpha_{ij}$. To obtain the relations we adjoin the codimension-two strata, where two circles are nested. If these occur in two different planes we have commutation relations $\alpha_{ij}\alpha_{k\ell}=\alpha_{k\ell}\alpha_{ij}$.  If the two occurrences of nested circles occur in the same plane we have either the second or the fifth configuration in the preceding large diagram.  The fifth configuration gives another commutation relation $\alpha_{ik}\alpha_{jk}=\alpha_{jk}\alpha_{ik}$.  The second configuration gives a relation $\alpha_{ij}\alpha_{ik}\alpha_{jk}=\alpha_{jk}\alpha_{ik}\alpha_{ij}$.  \qed

\medskip

This argument does not immediately extend to the groups $PUW_n$ since the space $\PUW^0_n$ corresponding to $\PUR^0_n$ is not simply-connected.  Its fundamental group is the pure braid group $PB_n$, so in principle it should be possible to extend a presentation for $PB_n$ to a presentation for $PUW_n$ by adjoining the generators $\alpha_{ij}$ corresponding to the codimension-one strata as before, and then figuring out the relations that correspond to the codimension-two strata.

The argument in the preceding proof does however work to prove the following:

\begin{prop}
For the subspace $\UR^<_n $ of $\UR_n$ consisting of configurations of rings of unequal size, there is a presentation for $\pi_1\,\UR^<_n$ with generators the $\alpha_{ij}$'s with $i<j$ and with relations the same relations as in the preceding proposition, restricted to these generators.
\end{prop}

\vspace{-9pt}

\proof By ordering rings according to size we obtain an embedding $\UR^<_n\incl\PUR_n$. The argument is then similar to the one for $\PUR_n$.  A small adjustment is needed in showing the subspace of unnested configurations has the homotopy type of the space of ordered point configurations; this we leave to the reader.  \qed

\medskip

\begin{prop}
A presentation for the group $W_n$ is obtained from the earlier presentation for $UW_n$ by adding the generators $\tau_i$ for $1\le i\le n$ and the following relations:
%\smallskip
\begin{enumerate}
\item[]%{}
\hspace{-6pt} $\left[\tau_i,\tau_j\right]=1$ for $i\ne j$
\smallskip
\item[]%{}
\hspace{-6pt} $\left[\rho_i,\tau_j\right]=1$ \ and  \  $\left[\sigma_i,\tau_j \right]=1$ \  for $j\ne i, i+1$
\smallskip
\item[]%{}
\hspace{-6pt} $\tau_i^\e \sigma_i^\eta = \sigma_i^\eta \tau_{i+1}^\e$ \ and \ 
$\tau_{i+1}^\e \sigma_i^\eta = \sigma_i^\eta \tau_i^\e$  \  for $\e,\eta = \pm 1$
\smallskip
\item[]%{}
\hspace{-6pt} $\tau_i^\e \rho_i = \rho_i \tau_{i+1}^\e$ \ and  \ 
$\tau_{i+1}^\e\rho_i = \sigma_i^{-\e}\rho_i^{-1}\sigma_i^\e\tau_i^\e$ \  for $\e = \pm 1$
\smallskip
\item[]%{}
\hspace{-6pt} $\tau_i^\e\rho_i^{-1}=\sigma_i^{-\e}\rho_i\sigma_i^\e\tau_{i+1}^\e$ \  and   \ $\tau_{i+1}^\e \rho_i^{-1} = \rho_i^{-1} \tau_i^\e$  \  for $\e = \pm 1$
\end{enumerate}
\end{prop}

The relations in the last three lines are highly redundant.  For example, two of the eight relations in the third-to-last line imply the other six. 

\vspace{-3pt}
\proof   It is not difficult to verify that the relations listed in the statement hold. These relations guarantee that any product of $\rho_i$'s, $\sigma_i$'s, and $\tau_i$'s can be rearranged as a product $ut$ where $u$ is a product of $\rho_i$'s and $\sigma_i$'s and $t$ is a product of $\tau_i$'s.

To verify that the $\rho_i$'s, $\sigma_i$'s, and $\tau_i$'s generate $W_n$ note first that for a given $x\in W_n$ there exists a product $s$ of $\sigma_i$'s such that $sx$ is in the subgroup $PW_n$. As we saw at the end of the preceding section, $PW_n$ is a semidirect product of $PUW_n$ and the subgroup $\Z^n$ generated by the $\tau_i$'s.  Thus $sx=ut$ for some $u\in PUW_n$ and $t$ a product of $\tau_i$'s.  Since $u$ is in $PUW_n$ it is in $UW_n$ and can therefore be written as a product of $\rho_i$'s and $\sigma_i$'s since we know these generate $UW_n$.  This implies that $x=s^{-1}ut$ is a product of $\rho_i$'s, $\sigma_i$'s, and $\tau_i$'s, so these elements generate $W_n$.  

To prove that the relations listed (including those for $UW_n$) define $W_n$, it will suffice to show that a word $w$ in the generators that represents the trivial element of $W_n$ can be reduced to the trivial word by applying the relations.  To start, we can use the relations to rewrite $w$ in the form $ut$ where $u$ is a product of $\rho_i$'s and $\sigma_i$'s (thus $u\in UW_n$) and $t$ is a product of $\tau_i$'s.  Since $ut=1$ and the $\tau_i$'s do not permute the wickets, we see that $u$ in fact lies in $PUW_n$.  The relation $ut=1$ implies that $u=1$ and $t=1$ in view of the semidirect product structure on $PW_n$.  The relations for $UW_n$ then suffice to reduce $u$ to the trivial word, and the commutation relations among the $\tau_i$'s allow $t$ to be reduced to the trivial word since the relation $t=1$ holds in the group $\Z^n$.  \qed

\medskip

The same argument works also for $\vr_n$ and $R_n$ to prove:

\begin{prop}
Presentations for $R_n$ and $\vr_n$ are obtained from the presentation for $W_n$ by adding the relations $\sigma_i^2=1$ and $\tau_i^2=1$ for $R_n$, or just $\sigma_i^2=1$ for $\vr_n$. \qed

\end{prop}

\vspace{-3pt}
Note that the relations involving the $\tau_i$'s can be simplified when $\sigma_i=\sigma^{-1}_i$.

\section{Rigidifying Floppy Wickets and Rings.}
\label{rigidifying}

Generalizing the spaces $\A_n$ and $\LL_n$ there is a space $\AL_{m,n}$ of smoothly embedded configurations of $m$ arcs and $n$ loops in $\R^3_+$ which are unknotted and unlinked, hence are isotopic to a configuration in $\WR_{m,n}$.  We also require the loops to be disjoint from the $xy$-plane.  Thus $\AL_{m,0} = \A_m$, and $\AL_{0,n}$ is homeomorphic to $\LL_n$ by the same argument that showed that $\WR_{0.n}$ is homeomorphic to $\RR_n$.

\begin{thm}
\label{3.1} 
The inclusion $\WR_{m,n}\incl \AL_{m,n}$ is a homotopy equivalence. \end{thm}

\vspace{-3pt}
Note that Theorems 1 and 2 in the Introduction follow directly as corollaries of Theorem~\ref{3.1}.

\vspace{-4pt}
\proof The space $ \WR_{m,n} $ is a smooth manifold and hence a CW complex, and $ \AL_{m,n} $ has the homotopy type of a CW complex, so it will suffice to show that the relative homotopy groups $ \pi_k(\AL_{m,n},\WR_{m,n}) $ vanish. As noted in the remarks following the proof of Theorem \ref{1.1}, the inclusion $\WR^c_{m,n}\incl\WR_{m,n}$ is a homotopy equivalence for each $c>0$, so it will in fact suffice to deform a given map $ f\cln (D^k,\bdy D^k) \to (\AL_{m,n},\WR^c_{m,n}) $ through such maps to a map $(D^k,\bdy D^k) \to (\WR_{m,n},\WR^c_{m,n})$, for any convenient choice of $c>0$.

Denote the family of arc and loop systems $f(t)$ by $ A_t $. We will be interested in systems $ S_t $ consisting of finitely many disjoint smooth disks and spheres embedded in $ \R^3_+ - A_t $ with $S_t\cap \bdy\R^3_+ = \bdy S_t$, such that each component of $ \R^3_+ - S_t $ contains at most one component of $ A_t $. We call such systems {\it separating systems}.  We assume that for each component of $ S_t $ there is a connected open set in the parameter domain $ D^k $ such that the component of $S_t$ varies only by isotopy as $ t $ ranges over this open set, and outside the open set the component is deleted from $ S_t $. If we choose the constant $c$ in $\WR^c_{m,n}$ to be less than $\frac{1}{2n}$ then for $ t \in \bdy D^k $ we can choose $S_t$ to consist of at least one round hemisphere or sphere in the interior of the microcosm of each wicket or ring of $A_t$, lying outside the wicket or ring, concentric with it, and disjoint from all other wickets and rings of $A_t$. By the remarks following the proof of Theorem~\ref{1.1},
we can assume that microcosms are disjoint from larger circles (and wickets).  This prescription for $S_t$ gives a separating system since each hemisphere or sphere chosen separates the corresponding wicket or ring from all other wickets or rings of equal or larger radius.  For nearby $t$ in $ \bdy D^k$ the hemispheres and spheres of the same radii remain a separating system, so we obtain in this way a family of separating systems $S_t$ consisting of round hemispheres and spheres for all $t$ in $\bdy D^k$.

There will be three main steps in the proof:

\begin{list}{}{\setlength{\leftmargin}{22pt}\setlength{\labelwidth}{16pt}\setlength{\labelsep}{5pt}\setlength{\itemsep}{3pt}}

\item[(1)] Extend the family of round separating systems $ S_t $ over $ \bdy D^k $ to smooth separating systems $ S_t $ for $ t \in D^k$.

\item[(2)] Deform these smooth separating systems to be round spheres and hemispheres over all of $ D^k $.

\item[(3)] Deform $ A_t $ so that it consists of round wickets and rings over all of $ D^k $. 
\end{list}
%\smallskip\noindent
At each step the family $A_t$ over $\bdy D^k$ will be unchanged.

\smallskip\noindent
{\bf Step 1: Extending over the disk.} There is a fibration $ \Diff(\R^3_+) \to \AL_{m,n} $ that sends a diffeomorphism to the image of a standard configuration of arcs and circles under the diffeomorphism. Using the lifting property of this fibration, we can choose a separating system for one parameter value $ t \in D^k $ and extend this to a family of separating systems $ \Sigma_t $ for $ A_t $ that varies only by isotopy as $ t $ ranges over all of $ D^k $. For $t \in \bdy D^k $ we then have two families of separating systems $ S_t $ and $ \Sigma_t $, and it will suffice to construct a family $ S_{tu} $, $ (t,u) \in \bdy D^k \times I $, which for each $ u $ is a separating system for $ A_t $, such that $S_{t0} = S_t $ and $ S_{t1} = \Sigma_t $.  We can then place this family $S_{tu}$ in a collar neighborhood of $\bdy D^k$ in $D^k$, after first deforming the family $A_t$ to be constant on each radial segment in this collar.

First thicken $ \Sigma_t $ to a family $ \Sigma_t \times [-1,1] $ of parallel separating systems for $ A_t $. Sard's theorem implies that for each $ t \in \bdy D^k $ there is a slice $ \Sigma_t \times \{ s \} $ in this thickening that is transverse to $ S_t $. This slice will remain transverse to $ S_t $ for all nearby $ t $ as well. By a compactness argument this means we can choose a finite cover of $ \bdy D^k $ by open sets $ U_i $ so that $ S_t $ is transverse to a slice $ \Sigma_i = \Sigma_i(t) $ for all $ t \in U_i $. 

For a fixed $t\in U_i$ consider the standard procedure for surgering $ S_t $ to make it disjoint from $ \Sigma_i $. The procedure starts with a component of $ S_t \cap \Sigma_i $, either a circle or an arc, that cuts off a disk $ D $ in $ \Sigma_i $ that contains no other components of $ S_t \cap \Sigma_i $. Using $ D $ we then surger $ S_t $ to eliminate the given component of $ S_t \cap \Sigma_i $. The process is then repeated until all components have been eliminated. Note that each surgery produces a system of disks and spheres that still separates $ \R^3_+ - A_t $ into components each containing at most one component of $ A_t $. 

A convenient way to specify the order in which to perform the sequence of surgeries is to imagine the surgeries as taking place during a time interval, and then surgering an arc or circle at the time given by the area of the disk it cuts off in $ \Sigma_i $, normalized by dividing by the area of $ \Sigma_i $ itself. The only ambiguity inherent in this prescription occurs if one is surgering the last remaining arc and this arc splits $ \Sigma_i $ into two disks of equal area. Then one would have to make an arbitrary choice of one of these disks as the surgery disk.

We will refine this procedure so that it works more smoothly in our situation. Thicken $ S_t $ to a family $ S_t \times [-1,1] $ of nearby parallel systems, all still transverse to $ \Sigma_i $ for $ t \in U_i $. Call this family of parallel systems $ {\bf S}_t $. For $ t \in U_i $, with $ i $ fixed for the moment, we perform surgery on $ {\bf S}_t $ by gradually cutting through it in a neighborhood of $ \Sigma_i $, as shown in Figure~4. 
\begin{figure}[htp]
\begin{center}
\includegraphics[width=2in]{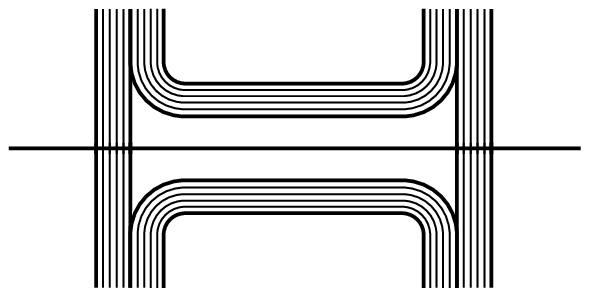}
\caption{ }
\label{ fig5}
\end{center}
\end{figure}
Thus we are producing a family $ {\bf S}_{tu} $ for $ u \in [0,1] $, where again we use the areas of the surgery disks in $ \Sigma_i $ to tell when to perform the surgeries. Notice that $ {\bf S}_{tu} $ is allowed to contain finitely many pairs of spheres or disks that touch along a common subsurface at the instant when these spheres or disks are being surgered. To specify the surgeries more completely we choose a small neighborhood $ \Sigma_i \times (-\varepsilon_i,\varepsilon_i) $ of $ \Sigma_i $ in $ \Sigma_t \times [-1,1] $, which we rewrite as $ \Sigma_i \times \R $, and we let the surgery on a component surface of $ {\bf S}_{tu} $ produce two parallel copies of the surgery disk in the slices $ \Sigma_i \times \{\pm 1/u\} $ of $ \Sigma_i \times \R $. Observe that this prescription for constructing $ \bf S_{tu} $ avoids the ambiguity in choosing one of the two equal-area surgery disks mentioned earlier since we can now surger using both these disks simultaneously. 

To convert the thickened family $ {\bf S}_{tu} $ back into an ordinary family $ S_{tu} $ consisting of finitely many disks and spheres for each $ (t,u) $ we replace each family of parallel disks or spheres in $ {\bf S}_{tu} $ of nonzero thickness by the central disk or sphere in this family. Thus this central disk or sphere belongs to $ S_{tu} $ for an open set of values of $ (t,u) $. 

As $ t $ varies over $ U_i $ we now have a family $ S_{tu} $, depending on $ i $. To combine these families for different values of $ i $, letting $ t $ range over all of $ \bdy D^k $ rather than just over $ U_i $, we proceed in the following way. For each $ i $ choose a continuous function $ \varphi_i \cln U_i \to [0,1] $ that takes the value $ 1 $ near $ \bdy U_i $ and the value $ 0 $ on an open set $ V_i $ inside $ U_i $ such that the different $ V_i $'s still cover $ D^k $. Then construct $ S_{tu} $ by delaying the time when each surgery along $ \Sigma_i $ is performed by the value $ \varphi(t) $. We may assume all the systems $ \Sigma_i $ are disjoint for fixed $t$ and varying $i$ with $t\in U_i$, and the thickenings $ \Sigma_i \times (-\varepsilon_i,\varepsilon_i) $ are disjoint as well, so the surgeries along different $ \Sigma_i $'s are completely independent of each other.

We have constructed the family $ S_{tu}$ for $ (t,u) \in \bdy D^k \times [0,1] $ such that all the curves of $ S_t \cap \Sigma_i $ are surgered away as $u $ goes from $ 0 $ to $ 1/2 $ for $ t \in V_i $. We can then adjoin $ \Sigma_i $ to $ S_{tu} $ for $ (t,u) \in V_i \times (1/2,1) $, deleting the surgered disks and spheres of $ S_{tu} $ for $ u \ge 3/4 $. We may assume all the thickenings $ \Sigma_i \times (-\varepsilon_i,\varepsilon_i) $ are disjoint from the original separating system $ \Sigma_t $. Then we adjoin $ \Sigma_t $ to $ S_{tu} $ for $ u > 3/4 $, so that for $ u = 1 $ only $ \Sigma_t $ remains in $ S_{tu} $.  This finishes Step~1.

\smallskip\noindent
{\bf Step 2: Rounding smooth disk and sphere systems.} We will use the following result:

\begin{lem}
The space of systems of finitely many disjoint smooth disks and spheres in $\R^3_+$, where the disks have their boundaries in $\bdy\R^3_+$, deformation retracts onto the subspace of round disks and spheres.
\end{lem}

\vspace{-9pt}
\proof We show the relative homotopy groups are zero, which is all we need for the application of the lemma.  Thus we are given a family $S_t$, $t\in D^k$, of disjoint smooth disks and spheres that we wish to isotope to round disks and spheres, staying fixed over $\bdy D^k$ where $S_t$ is assumed to already consist of round disks and spheres. We can assume in fact that $S_t$ consists of round disks and spheres for $t$ in a neighborhood of $\bdy D^k$.  

First we show how to round the spheres of $S_t$ by an inductive procedure, starting with the outermost spheres.  We construct families of embeddings of $D^3$ in $\R^3_+$ with images bounded by the outermost spheres, such that near $\bdy D^k$ these embeddings are rescaled isometric embeddings. This can be done by first applying isotopy extension to construct families of embeddings without the condition near $\bdy D^k$, then deforming these embeddings to achieve this extra condition using the fact that the inclusion of $O(3)$ into $\Diff(D^3)$ is a homotopy equivalence, which is a consequence of the Smale conjecture that $ \Diff(D^3 \rel\bdy D^3) $ is contractible, proved in [H1].  We can also arrange that the embeddings are rescaled isometric embeddings near the center of $D^3$, just by differentiability.  By restricting these embedding to smaller and smaller concentric spheres in $D^3$ we can isotope the outermost spheres to be round over all of $D^k$, damping the isotopy down to the identity near $\bdy D^k$. The non-outermost spheres are dragged along in this process. Having rounded the outermost spheres in $S_t$, we do a similar construction for the next-outermost spheres, and so on.  

To make the disks round we first make all their boundary circles round following the same plan as for spheres, using Smale's theorem that $ \Diff(D^2 \rel \bdy D^2) $ is contractible.  The rounding of the boundary circles can be done by a deformation of the family $S_t$ supported in a neighborhood of $\bdy\R^3_+$.  Having the boundary circles round, we then deform the disks themselves to the round hemispherical disks spanning the round boundary circles. This is possible since the fibration obtained by restricting the disks to their boundaries has contractible fiber, the space of smooth disk systems in $ \R^3_+ $ with given boundary circles. For a single disk this is one of the equivalent forms of the Smale conjecture, and for systems of disks it follows by induction. When we perform these isotopies of the disks of $ S_t $, the spheres of $S_t$ are to be dragged along, so the proper way to proceed is first to make all the disks round, then make the spheres round by the procedure described earlier. 
\qed
\medskip

Now we return to Step 2 of the proof.  For each $ t_0 \in D^k $ the components of $ S_t $ vary only by isotopy as $t$ varies over some neighborhood of $ t_0 $. Choose a finite number of these neighborhoods that cover $ D^k $, then triangulate $ D^k $ so that each $k$-simplex of the triangulation lies in one of these neighborhoods. Over each such $k$-simplex we then have the associated set of disks and spheres of $ S_t $ that vary only by isotopy. Over a face of the simplex we have the union of the sets of disks for the various $ k $-simplices that contain the face. Let us change notation slightly and call these systems of surfaces $S_t $. (They are subsets of the systems $ S_t $ constructed in Step~1.) 

Suppose inductively that we have isotoped the disks and spheres of $ S_t $ to be round for $ t $ in the $ i $-skeleton of the triangulation of $ D^k $, without changing anything over $ \bdy D^k $ where the systems $ S_t $ and $ A_t $ are already round. The possibility $i=-1$ is allowed, which will give the start of the induction. For the induction step we apply the lemma to extend the rounding isotopy of $S_t$ over each $(i+1)$-simplex in the interior of $D^k$ in turn.  The arcs and circles of $A_t$ are carried along during this deformation of $S_t$, by isotopy extension. This completes Step~2.

\smallskip\noindent
{\bf Step 3: Rounding smooth arc and circle systems.}  Having the components of $ S_t $ round over all of $ D^k $, we can round the components of $ A_t $ by an inductive procedure as in Step~2. Over a simplex $\sigma$ of the triangulation of $ D^k $ we look at a complementary region $C_t$ of $ S_t $. This contains at most one component of $ A_t $, and we need only look at the case when there is exactly one component, say $ \alpha_t $. Consider first the case that $\alpha_t$ is an arc. The region $C_t$ is bounded by the plane $ \R^2 $ together with some round disks and spheres of $ S_t $ that can vary by isotopy. Let $\widehat C_t$ be obtained from $C_t$ by filling in the boundary spheres with balls.  We can then think of $\widehat C_t$ as a region in the upper half-space model of hyperbolic $3$-space bounded by geodesic planes. There is always a unique round arc $ \alpha'_t $ in $\widehat C_t$ having the same endpoints as $ \alpha_t $. This means that the space of round arcs in $\widehat C_t$ is the same as the space of pairs of endpoints of smooth arcs. The map sending each unknotted smooth arc to its endpoints is a fibration, and it is a homotopy equivalence since its fiber, the space of unknotted arcs with fixed endpoints, is contractible, by another equivalent form of the Smale conjecture. Since the fibration is a homotopy equivalence, this implies that we can deform the arcs $ \alpha_t $ to round arcs over the simplex $ \sigma $, staying fixed over the boundary of $\sigma$ where they are already round.  We can drag the balls of $\widehat C_t - C_t$ and everything inside them along during the isotopy that rounds $\alpha_t$. This could destroy the roundness of these balls, but this problem can be avoided by first shrinking the balls sufficiently small so that they can stay round during the isotopy.

The other case is that $\alpha_t$ is a circle. There are then two subcases depending on whether $C_t$ is of the same type as in the preceding case or $C_t$ is a ball with smaller disjoint sub-balls removed. In the first subcase the space of round circles in $\widehat C_t$ has the homotopy type of $\R P^2$ since such circles bound unique geodesic disks in $\widehat C_t$ and the space of such disks has this homotopy type. The space of smooth unknotted circles in $\widehat C_t$ also has the homotopy type of $\R P^2$ by the Smale conjecture, so we can deform the circles $\alpha_t$ to be round over $\sigma$ as before, after first shrinking the balls of $\widehat C_t - C_t$. The other subcase, that $C_t$ is a ball with sub-balls removed, is done in the same way, using the fact that the space of round circles in a ball has the same homotopy type as the space of smooth circles, namely $\R P^2$ again.    

This finishes the proof of the theorem.  \qed
\smallskip

\p{Further injectivity results.} We observed at the beginning of Section~\ref{presentations} that $W_n$ contains two copies of the braid group $B_n$, one generated by the $\rho_i$'s and the other generated by the $\sigma_i$'s.  Under the projection $W_n\to R_n$ the copy of $B_n$ generated by the $\sigma_i$'s becomes a subgroup $\Sigma_n \subset R_n$, and we can now see that the other copy of $B_n$ remains unchanged:

\begin{prop}
\label{3.2} 
The map $\sigma : B_n\to R_n$ sending the standard generators of the braid group to the elements $\rho_i$ is injective.
\end{prop}

\vspace{-9pt}
\proof It suffices to show $\sigma$ is injective on the `pure' versions of these groups, the kernels of the natural maps to $\Sigma_n$.  The pure braid group fits into a well-known split short exact sequence
$$
0\to F_{n-1}\to PB_n \to PB_{n-1}\to 0
$$
where $F_{n-1}$ is the free group on $n-1$ generators and the map $PB_n\to PB_{n-1}$ is obtained by ignoring the last strand of a pure braid.   This short exact sequence maps to a similar split short exact sequence
$$
0\to K_n\to PR_n\to PR_{n-1}\to 0
$$
which is part of the long exact sequence of homotopy groups associated to the fibration which sends an ordered $n$-tuple of smooth circles forming the trivial link to the ordered $(n-1)$-tuple obtained by ignoring the last circle.  The kernel $K_n$ is $\pi_1$ of the fiber, the subspace of $\LL_n$ consisting of configurations with $n-1$ of the circles in a fixed position and the last circle varying.  It suffices by induction on $n$ to show that the map of kernels $F_{n-1}\to K_n$ is injective.  We do this by constructing a homomorphism $K_n\to F_{n-1}$ such that the composition $F_{n-1}\to K_n\to F_{n-1}$ is the identity.  

The homomorphism $K_n\to F_{n-1}$ is obtained by choosing a point in the $n$th circle and taking the path it traces out in the complement of the other $n-1$ circles under a loop in the fiber.  This path may not be a loop, but it can be completed to a loop by adjoining an arc in the $n$th circle.  Since the circles are unlinked, the choice of this arc does not affect the resulting element of $F_{n-1}$, the fundamental group of the complement of the first $n-1$ circles.  This construction gives a homomorphim $K_n \to F_{n-1}$ such that precomposing with $F_{n-1}\to K_n$ is obviously the identity. \qed

\medskip

The kernel $K_n$ is the product $KU_n\times \Z$ for $KU_n$ the kernel of the projection ${\pur}_n\to {\pur}_{n-1}$.  It is shown in [P] that $KU_n$ is not finitely presented for $n\ge 3$, although it is finitely generated, with the generators one might expect, $\alpha_{ni}$ and $\alpha_{in}$.  The lack of finite presentability probably means that these kernels do not have nice geometric interpretations in terms of configuration spaces of circles.

\begin{prop}
\label{3.3}
The map ${\ur}^<_n\to {\ur}_n$ is injective.
\end{prop}

\vspace{-9pt}
\proof  This is similar to the preceding proof. The map $\UR^<_n\to \UR^<_{n-1}$ that ignores the smallest ring is a quasifibration, as in Section~1, using the canonical shrinking to first make the smallest ring point-sized.  The fundamental group of the fiber is $F_{n-1}$ so we get a split short exact sequence
$$
0\to F_{n-1}\to {\ur}^<_n\to {\ur}^<_{n-1}\to 0
$$
which maps to the split short exact sequence
$$
0\to K_n \to PR_n\to PR_{n-1}\to 0
$$
from the preceding proof. The rest of the argument is the same. \qed

\medskip

\section{Asphericity.}
\label{asphericity}

As a warm-up to proving Theorem 3, which states that the spaces $\W_n$ and $\UW_n$ are aspherical, let us recall a standard sort of argument for showing that the map $W_n\to B_{2n}$ induced by the map $\A_n\to\C_{2n}$ sending a configuration of arcs to the configuration of its endpoints is injective.   We can view $\A_n$ as the space of configurations of $n$ disjoint smooth unknotted, unlinked arcs in a ball $D^3$ with endpoints in a hemisphere $D^2_-$ of $\bdy D^3$.  By restricting diffeomorphisms of $D^3$ fixing the other hemisphere $D^2_+$ to the standard configuration $A$ of $n$ arcs we obtain a fibration 
$$
\Diff(D^3, A\rel D^2_+) \to \Diff(D^3\rel D^2_+) \to \A_n \eqno(1)
$$
where $\Diff(X,Y \rel Z)$ denotes the space of diffeomorphisms of a manifold $X$ that leave a submanifold $Y$ setwise invariant and fix a submanifold $Z$ pointwise.  Restricting everything to $D^2_-$ gives a map from this fibration to the fibration
$$
\Diff(D^2_-,\bdy A\rel \bdy D^2_-) \to \Diff(D^2_- \rel\bdy D^2_-) \to \C_{2n} \eqno(2)
$$
In each fibration the projection map to the basespace is nullhomotopic by shrinking the support of diffeomorphisms to a smaller ball or disk disjoint from $A$. Thus the associated long exact sequences of homotopy groups break up into short exact sequences. Since $\pi_0\Diff(D^2_-\rel\bdy D^2_-)=0$ and $\pi_0\Diff(D^3\rel D^2_+)=0$ (the latter by Cerf's theorem), we obtain isomorphisms $\A_n\approx\pi_0\Diff(D^3, A\rel D^2_+)$ and $B_{2n}\approx\pi_0\Diff(D^2_-,\bdy A\rel \bdy D^2_-)$. The problem is thus reformulated as showing injectivity of the map 
$$
\pi_0\Diff(D^3, A\rel D^2_+)\to \pi_0\Diff(D^2_-,\bdy A\rel \bdy D^2_-)
$$ 
This map is induced by the restriction map from the fiber of the first fibration above to the fiber of the second fibration.  This restriction map is itself a fibration 
$$
\Diff(D^3,A\rel\bdy D^3)\to \Diff(D^3, A\rel D^2_+) \to \Diff(D^2_-,\bdy A\rel \bdy D^2_-) \eqno(3)
$$
so it suffices to show that $\pi_0$ of the fiber of this fibration is trivial. Note first that a diffeomorphism $f$ in $\Diff(D^3,A\rel\bdy D^3)$ can be isotoped to be the identity on $A$, and $f$ cannot twist the normal bundles of the arcs of $A$, as one can see by looking at the induced map on $\pi_1(D^3-A)$. Then $f$ can be isotoped rel $A\cup \bdy D^3$ to be the identity in a neighborhood of $A$, so $f$ can be regarded as a diffeomorphism of a handlebody fixing the boundary of the handlebody. The space of such diffeomorphisms is path-connected since any two spanning disks in a handlebody are isotopic rel boundary, and similarly for collections of disjoint spanning disks, so diffeomorphisms of a handlebody rel boundary can be isotoped rel boundary to have support in a ball, and then by Cerf's theorem they can be isotoped to the identity.   (With a little more work the use of Cerf's theorem in this argument could be avoided by factoring out the image of $\pi_0\Diff(D^3\rel\bdy D^3)$ in the various groups.)

\medskip
Now we prove Theorem 3 by refining this argument to reduce asphericity of $\W_n$ to asphericity of $\C_{2n}$.  

\vspace{-3pt}
\proof Since $\W_n$ is homotopy equivalent to $\A_n$, we can obtain the result for $\W_n$ by showing that $\A_n$ is aspherical.  The total space in the fibration (2) above is contractible by a theorem of Smale. The total space in the fibration (2) is also contractible, as one can see from the fibration
$$
\Diff(D^3\rel \bdy D^3)\to \Diff(D^3\rel D^2_+) \to \Diff(D^2_-\rel \bdy D^2_-)
$$
where the base is contractible by Smale's theorem and the fiber is contractible by the Smale conjecture [H1].  The fiber of the fibration (3) is also contractible by the following argument. Restricting diffeomorphisms in $ \Diff(D^3,A\rel\bdy D^3) $ to normal bundles of the $ n $ arcs gives another fibration whose base space is homotopy equivalent to the space of automorphisms of the normal bundles of these arcs that are the identity at the endpoints of the arcs. For each arc this is the loopspace of $ SO(2) $, which has contractible components. Components other than the identity component can be ignored since diffeomorphisms in $ (D^3,A\rel\bdy D^3) $ cannot twist the normal bundles nontrivially, as we saw earlier. Thus from this fibration we can replace $ (D^3,A\rel\bdy D^3) $ by the subspace of diffeomorphisms that are the identity on a neighborhood of the arcs. This can be identified with group of diffeomorphisms of a handlebody fixing its boundary. This diffeomorphism group is path-connected as we observed before, and it has contractible path-components by [H2]. (The key point is that the space of spanning disks with fixed boundary is contractible.)

Thus for $i\ge 2$ we have isomorphisms 
$$
\pi_i\A_n \approx \pi_{i-1}\Diff(D^3, A\rel D^2_+)\approx \pi_{i-1}\Diff(D^2_-,\bdy A\rel \bdy D^2_-)\approx \pi_i\C_{2n}
$$
so asphericity of $\A_n$ is reduced to asphericity of $\C_{2n}$, which is well-known.  

For the case of $\UW_n$ we can pass to the covering space $\PUW_n$ obtained by ordering the wickets, and then use the quasifibration $\PUW_n\to \PW_n \to T^n$ from Section~1, where $T^n$ is the $n$-torus.  The associated long exact sequence of homotopy groups shows that $\PUW_n$ is aspherical since $\PW_n$ and $T^n$ are aspherical.   \qed

\medskip

\section{Wickets and Rings in a Sphere.}
\label{sphere}

Instead of wickets in upper halfspace one can consider wickets inside a sphere, circular arcs in the interior of the sphere that meet the sphere orthogonally at their endpoints.  Configurations of $n$ disjoint wickets of this type form a spherical wicket space $\SW_n$.  An equivalent space is the space of configurations of $n$ disjoint line segments in a ball that meet the boundary sphere in their endpoints. The equivalence between the two definitions can be seen by considering two of the models for hyperbolic $3$-space, the standard ball model and the projective model.  In the ball model the geodesics are circular arcs orthogonal to the boundary sphere, while in the projective model they are line segments in the ball with endpoints on the boundary sphere.  The disjointness condition is preserved in going from one model to the other since intersecting geodesics lie in a common hyperbolic plane in both cases.  

The space $\SW_n$ is slightly smaller than the space of all configurations of $n$ disjoint geodesics in hyperbolic $3$-space since geodesics do not include their endpoints in the boundary sphere, so two disjoint geodesics could share a common endpoint on the boundary sphere.  The inclusion of $\SW_n$ into this slightly larger space is a homotopy equivalence, however, as one can see easily in the projective model by shrinking the ball by a small amount for each configuration (without shrinking the configuration itself).  For example, the ball can be shrunk by one-half of the minimum of the numbers $d_i$, where $d_i$ is the maximum distance from points on the $i$th line segment of a given configuration to the boundary of the ball. Note that this is essentially the same as the canonical shrinking process considered in Section~1.

Comparing the ball model of hyperbolic $3$-space with the upper halfspace model, we see that $\W_n$ can be regarded as the subspace of $\SW_n$ consisting of configurations disjoint from a point $\infty$ in the boundary sphere.  The configurations in $\SW_n$ that contain a line to $\infty$ form a codimension $2$ submanifold.   In terms of the upper halfspace model, this submanifold is the space of configurations of $n-1$ disjoint wickets and one vertical line disjoint from the wickets. 
This submanifold is connected, by the same argument with canonical shrinking used to show that $W_n$ is connected.  From transversality it follows that the inclusion $\W_n \incl \SW_n$ induces a surjection on $\pi_1$ with kernel generated by a small loop linking the codimension $2$ submanifold.  This loop can be represented by taking the standard configuration of $n$ wickets in the $xz$-plane and dragging the left endpoint of the first wicket around a large circle enclosing all the other wickets.  It would not be hard to write this loop as a word in the generators $\rho_i$, $\sigma_i$ and $\tau_i$.  Thus $\pi_1\SW_n$ has a presentation obtained from the presentation for $W_n$ by adding one extra relation.

There is an analogous space $\SA_n$ of configurations of $n$ disjoint smooth arcs in a ball with endpoints on the boundary sphere, all these arcs being unknotted and unlinked. 

\begin{prop}\label{5.1}
The inclusion $\SW_n\incl\SA_n$ is a homotopy equivalence.
\end{prop}

\vspace{-9pt}
\proof This can be reduced to the corresponding result for $\W_n \incl \A_n$ by considering some fibrations. Let $\SW_N^*$ be the space of configurations consisting of $n$ disjoint wickets in a ball together with a basepoint in the boundary sphere disjoint from the wickets.  Projecting such a configuration onto either the wickets or the basepoint gives two fibrations
$$
F\to \SW_n^*\to \SW_n \qquad\qquad
\W_n\to\SW_n^*\to S^2
$$
Here the fiber $F$ in the first fibration is just $S^2$ with $2n$ points deleted, the endpoints of a configuration of $n$ wickets.  The homotopy lifting property in the first fibration follows by extending isotopies of configurations of wickets to ambient isotopies then restricting these to the basepoint. The second fibration is actually a fiber bundle since the basepoints in a neighborhood of a given basepoint can be obtained via a continuous family of rotations of $S^2$ applied to the given basepoint, and then these rotations can be applied to configurations of wickets.

Similarly there are fibrations
$$
F\to \SA_n^*\to \SA_n \qquad\qquad
\A_n\to\SA_n^*\to S^2
$$
The fiber $F$ is the same as before.  There are natural maps from the first two fibrations to the second two fibrations.  Applying the five lemma to the induced maps of long exact sequences of homotopy groups, we see that $\W_n \incl \A_n$ being a homotopy equivalence implies first that this is true also for $\SW_n^*\incl\SA_n^*$ and then also for $\SW_n\incl\SA_n$. \qed

\medskip
Similar things can be done for rings as well as wickets.  Let $\SR_n$ be the space of configurations of $n$ disjoint pairwise unlinked circles in $S^3$, and let $\SL_n$ be the corresponding analog of $\LL_n$, the space of smooth $n$-component trivial links in $S^3$.  

\begin{prop}
The inclusion $\SR_n\incl\SL_n$ is a homotopy equivalence.
\end{prop}
\vspace{-9pt}
\proof This follows the line of argument in the preceding proof by comparing fibrations, using the space $\SR_n^*$ of configurations of circles in $S^3$ with a disjoint basepoint, and its smooth analog $\SL_n^*$. \qed

\medskip
One can also obtain a presentation for $\pi_1\SR_n$ from a presentation for $R_n$ by adding the same relation as was added to get a presentation for $\pi_1\SW_n$.  The justification is the same as before, by using stereographic projection to identify $\RR_n$ with the complement of the codimension $2$ submanifold of $\SR_n$ consisting of configurations passing through a given point in $S^3$.

\section{Remarks on Dimension}
\label{dimension}

It is a classical fact that the general position argument for finding a presentation for $B_n$ can be refined to build a finite CW complex $K(B_n,1)$ having a single $0$-cell, a $1$-cell for each standard generator $\sigma_i$, and a $2$-cell for each of the standard relations. 
The cells are dual to the strata of the stratification of $\C_n$ according to coincidences of the $x$-coordinates. Thus the $0$-cell corresponds to the unique stratum of maximum dimension consisting of configurations with distinct $x$-coordinates, the $1$-cells to the strata of codimension one where exactly two points in a configuration have the same $x$-coordinate, and so on.  The same procedure works also for $\UW_n$ to give a finite CW complex $K(UW_n,1)$. The dimension of this complex is $n-1$,  just as for $B_n$. For $B_n$ there is a single cell in the top dimension, corresponding to the stratum of configurations with all $n$ points on one vertical line, but for $\UW_n$ there are a number of different strata consisting of configurations of wickets all lying in one plane, so there are a number of top-dimensional cells.  There cannot exist a $K(UW_n,1)$ of dimension less than $n-1$ since $UW_n$ has a subgroup $\Z^{n-1}$ generated by the elements $\alpha_{in}$ for $i<n$.

For $W_n$ the minimum dimension of a $K(W_n,1)$ is $2n-1$. There is a $K(W_n,1)$ of this dimension since $W_n$ is a subgroup of $B_{2n}$, and there cannot be one of lower dimension since $W_n$ contains a subgroup $\Z^{2n-1}$, generated by the $\Z^{n-1}$ above and the $\tau_i$'s. It seems likely that $W_n$ should have a finite CW complex $K(W_n,1)$ of minimum dimension, perhaps constructible by extending the general-position constructions referred to above. 

For $R_n$ the virtual cohomological dimension is known to be $n-1$ by [C], where a $K(\pi,1)$ which is a finite CW complex of dimension $n-1$ was constructed for the finite-index subgroup ${\pur}_n$. This $K(\pi,1)$ can be described as the space of basepointed graphs consisting of $n$ circles touching in a tree-like pattern, forming a cactus-shaped object. The dimension $n-1$ cannot be reduced since ${\pur}_n$ again contains a subgroup $\Z^{n-1}$ generated by the elements $\alpha_{in}$.

\section{References}

\parindent 45pt

\item{[BWC]} J.\,C.\,\,Baez, D.\,K.\,\,Wise, A.\,S.\,\,Crans, Exotic statistics for strings in 4d BF theory. {\em Adv.\ Theor.\ Math.\ Phys.}\ 11 (2007), 707--749. arXiv:gr-qc/0603085.

\item{[BMMM]} N.\,\,Brady, J.\,\,McCammond, J.\,\,Meier, and A.\,\,Miller. The pure symmetric automorphisms of a free group form a duality group. {\em J.\ Algebra} 246 (2001), 881--896.

\item{[BL]} A.\,\,Brownstein and R.\,\,Lee. Cohomology of the group of motions of n strings in 3-space. {\em Contemp.\ Math.}\ 150 (1993), 51--61. 

\item{[C]} D.\,J.\,\,Collins, Cohomological dimension and symmetric automorphisms of a free group.  {\em Comment.\ Math.\ Helv.}\  64  (1989), 44--61. 

\item{[CPVW]} F.\,R.\,\,Cohen, J.\,\,Pakianathan, V.\,\,Vershinin, and J.\,\,Wu, Basis-conjugating automorphisms of a free group and associated Lie algebras, arXiv:math.GR/0610946.

\item{[D]} D.\,M.\,\,Dahm, A Generalization of Braid Theory, Ph.D.~Thesis, Princeton University, 1962.

\item{[FRR]} R.\,\,Fenn, R.\,\,Rim\'anyi and C.\,\,Rourke, The braid-permutation group, {\em Topology} 36 (1997), 123--135.

\item{[FS]} M.\,\,Freedman and R.\,\,Skora, Strange Actions of Groups on Spheres, {\em J.\ Diff.\ Geom.}\ 25 (1987), 75--98.

\item{[G]} D.\,L.\,\,Goldsmith, The theory of motion groups. {\em Michigan Math.\ J.}\ 28 (1981), 3--17. 

\item{[H1]} A.\,\,Hatcher, A proof of the Smale conjecture, {\em Ann.\ of Math.}\ 117 (1983), 553--607.

\item{[H2]} A.\,\,Hatcher, Homeomorphisms of sufficiently large $P^2$-irreducible $3$-manifolds,  {\em Topology}  15  (1976), 343--347. Revised and updated version: Spaces of incompressible surfaces, arXiv:math.GT/9906074.

\item{[H3]} H.\,\,Hilden, Generators for two groups related to the braid group, {\em Pac.\ J.\ Math.}\ 59 (1975), 475--486.

\item{[JMM]} C.\,\,Jensen, J.\,\,McCammond, and J.\,\,Meier, The integral cohomology of the group of loops, {\em Geom.\ \& Top.}\ 10 (2006), 759--784.
 
\item{[Mc]} J.\,\,McCool, On basis-conjugating automorphisms of the free groups, {\em Can.\ J.\ Math.}\ 38 (1986), 1525--1529.

\item{[P]} A.\,\,Pettet, Finiteness properties for the kernel of pure motions of $n$ unlinked loops in $\R^3$, arXiv:math/0602148.

\item{[R]} R.\,L.\,\,Rubinsztein, On the group of motions of oriented, unlinked and unknotted circles in $\R^3$, I. Preprint, Uppsala University, 2002. Available online at http://www.math.uu.se/research/pub/Rubinsztein1.pdf

\item{[T1]} S.\,\,Tawn, A presentation for Hilden's subgroup of the braid group, {\em Math.\ Res.\ Lett.}\ 15 (2008), no.\ 6, 1277--1293.

\item{[T2]} S.\,\,Tawn, A presentation for the pure Hilden group, arXiv:0902.4840.

\item{[V]} V.\,V.\,\,Vershinin, On homological properties of singular braids, {\em Trans.\ A.M.S.}\ 350 (1998), 2431--2455.

\end{document}